\documentclass[a4paper,fleqn]{cas-sc}

\usepackage{algorithm}
\usepackage{algpseudocode}
\usepackage{enumitem}
\usepackage{xcolor} 
\usepackage{bm} 

\usepackage[numbers, sort&compress]{natbib}

\usepackage{graphicx} 
\usepackage{subcaption} 

\usepackage{float} 
\usepackage{tabularx}  
\usepackage{booktabs}  
\usepackage{multirow}  
\usepackage{makecell}

\def\tsc#1{\csdef{#1}{\textsc{\lowercase{#1}}\xspace}}
\tsc{WGM}
\tsc{QE}

\begin{document}
\let\WriteBookmarks\relax
\def\floatpagepagefraction{1}
\def\textpagefraction{. 001}
\let\printorcid\relax 

\shorttitle{DD-DeepONet}   

\shortauthors{Yang et al. }

\title[mode = title]{DD-DeepONet: Domain decomposition and DeepONet for solving partial differential equations in three application scenarios}

\author[1,2]{Bo Yang}

\author[2]{Xingquan Li}
\author[2]{Jie Zhao}


\author[1]{Ying Jiang}
\cormark[1]
\ead{jiangy32@mail.sysu.edu.cn} 

\address[1]{School of Computer Science and Engineering, Sun Yat-sen University, Guangzhou, 510006, Guangdong, China}
\address[2]{Pengcheng Laboratory, Shenzhen, 518052, Guangdong, China}

\cortext[1]{Corresponding author} 

\begin{abstract}
In certain practical engineering applications, there is an urgent need to perform repetitive solving of partial differential equations (PDEs) in a short period. This paper primarily considers three scenarios requiring extensive repetitive simulations. These three scenarios are categorized based on whether the geometry, boundary conditions(BCs), or parameters vary. We introduce the DD-DeepONet, a framework with strong scalability, whose core concept involves decomposing complex geometries into simple structures and vice versa. We primarily study complex geometries composed of rectangles and cuboids, which have numerous practical applications. Simultaneously, stretching transformations are applied to simple geometries to solve shape-dependent problems. This work solves several prototypical PDEs in three scenarios, including Laplace, Poission, N-S, and drift-diffusion equations, demonstrating DD-DeepONet’s computational potential. Experimental results demonstrate that DD-DeepONet reduces training difficulty, requires smaller datasets and VRAM per network, and accelerates solution acquisition.

\end{abstract}



\maketitle

\section{Introduction}
PDEs represent a class of differential equations that involve an unknown multivariate function along with its partial derivatives. These equations are paramount in modeling and describing an extensive array of physical phenomena, including but not limited to acoustic, heat, diffusion \cite{WOS:000513929200001}, electromagnetism \cite{Cai_2013}, fluid dynamics \cite{WOS:000501350300042}, mechanics \cite{WOS:000911026700001, WOS:000762477300009}. PDEs are pivotal in both scientific research and engineering applications. 

However, as is widely acknowledged, obtaining analytical solutions to PDEs can be exceptionally challenging. Most PDEs do not have straightforward analytical solutions and thus are predominantly solved using various numerical methods. The most used traditional numerical techniques including finite element method (FEM) \cite{WOS:000244987700052}, finite volume method (FVM) \cite{WOS:000408075700013}, finite difference method (FDM) \cite{INSPEC:811169}, smoothed particle hydrodynamics (SPH) \cite{ WOS:001343474100001}, lattice boltzmann method (LBM), molecular dynamics (MD) \cite{WOS:000993135400001}, dissipative particle dynamics (DPD) \cite{WOS:000399904500001}, spectral method \cite{shenjie_2011}, boundary element method (BEM) and so on. These methods provide approximate solutions that closely mimic the behavior described by the original PDEs. Each of these methods has its unique advantages and applicability depending on the specific requirements of the problem at hand. However, these numerical methods share a common challenge: their computational demands tend to increase exponentially with the increase in dimensionality and the refinement of the discretization of the computational region. This is especially true for 3D problems where achieving the desired level of accuracy may require significant computational time or resources, often necessitating the use of powerful supercomputers. The intricate balance between accuracy, computational time, and resource utilization continues to be a challenge in the numerical solution of PDEs. 

The concept of the domain decomposition method (DDM) was first proposed by the German mathematician H. A. Schwarz in 1870 \cite{Saad_DDM, AT_OBW_DDM}. However, it was not utilized in computation until the 1950s. With the advent and proliferation of parallel computing in the 1980s, DDM emerged as a solution to alleviate the computational limitations of traditional numerical methods. 

DDM employs a divide-and-conquer strategy. It breaks down a large problem into several smaller problems, solves each one independently, and then combines their solutions to form a comprehensive solution. This approach offers several advantages: 
\begin{enumerate}
    \item Reduced computational scale: By reducing the scale of each computational task, DDM significantly mitigates the limitations imposed by the capacity and speed of computers.
    \item Simplification of complex regions: complex computational domains can be decomposed into simpler or more regular subdomains, transforming complicated problems into simpler ones. 
    \item Utilization of fast algorithms: In each simpler subdomain, various fast algorithms can be applied, such as the Fast Fourier Transform (FFT) and spectral methods. 
    \item Flexible discretization: there is no need for a consistent mesh across the entire domain, different subdomains can employ different discretization methods. 
    \item Diverse equations in subdomains: completely different partial differential equations can be used in different subdomains, better reflecting the actual problem. 
    \item Highly parallel algorithms: DDM can be designed to be highly parallel, leveraging multiple computing resources to reduce computational time. 
\end{enumerate}
These characteristics make DDM particularly effective for large-scale scientific and engineering problems, enabling more efficient use of computational resources and facilitating the handling of complex multi-domain challenges. 

In recent years, advancements in technology, particularly the robust computational power of GPUs, have propelled the rapid development of artificial intelligence (AI) technologies, notably deep neural networks. Since the debut of OpenAI's ChatGPT, this domain has consistently been a focal point in both academic and industrial spheres. The awarding of the 2024 Nobel Prizes in physics and chemistry has brought unprecedented attention to the field of artificial intelligence, marking a significant milestone. 

As early as the 1990s, scholars were already exploring the mathematical foundations and methodologies for solving PDEs using neural networks. To date, researchers have employed a variety of neural networks for this purpose, including fully connected neural networks (FCNN), convolutional neural networks (CNN) \cite{726791}, residual neural networks (RNN), transformer \cite{10.5555/3295222.3295349}, and so on. Building on these foundational network architectures, numerous methods have been developed to address PDEs, such as the deep galerkin method \cite{WOS:000450907600062}, the deep ritz method \cite{WOS:000426076300001}, physics-informed neural networks (PINN) \cite{WOS:000453776000028}, operator learning \cite{WOS:001113803200001}, and so on. These developments highlight the intersection of deep learning and numerical computation, illustrating how AI is increasingly integral in advancing scientific computation and engineering applications. 

In engineering applications, solving PDEs often requires not just a one-time solution but repeated solutions based on practical needs. In this paper, we mainly consider three application scenarios, respectively:
\begin{enumerate}
    \item scenario 1 (S1): geometry fixed, BC and parameter varying. 
    \item scenario 2 (S2): geometry varying, BC, and parameter fixed.
    \item scenario 3 (S3): geometry, BC, and parameter all varying.
\end{enumerate}
For instance, the weather forecast is related to S1. In very large-scale integration (VLSI), parasitic resistance and capacitance extraction problems are associated with S2. For S3, topology and shape optimization, clinical prognostication, there is a significant demand for these analyses \cite{WOS:001372682600001, WOS:001256289700001}.

Using traditional numerical methods for solving PDEs can be computationally expensive and time-consuming, especially when dealing with 3D problems that involve a large number of collocation points and require extensive repeated calculations. Furthermore, a substantial amount of related data is generated after these simulations, which often remains underutilized, leading to a waste of data resources. 

To address S1, operator learning has been developed. Following the universal approximation theorem for functions proposed by Hornik K. and others in 1989 \cite{HORNIK1989359}, Chen T. and colleagues introduced the universal approximation theorem for nonlinear operators \cite{392253} in 1995. In recent years, several operator learning networks have been proposed, such as the nonlocal kernel network \cite{YOU2022111536}, Koopman neural operator \cite{XIONG2024113194}, learning operators with coupled attention \cite{JMLRv2321-1521}, and so on. Notably, the Fourier Neural Operator (FNO) \cite{li2021fourier} introduced by Li Z et al. and DeepONet \cite{WOS:000641834300001} by Lulu et al. have become particularly popular. Operator learning, capable of learning mappings between infinite-dimensional function spaces, effectively addresses the needs of S1 \cite{BAHMANI2025118113}. However, it remains ineffective for the S2 and S3, because of geometric encoding issues of collocation points and computation domain. Additionally, if the training data for the PDEs is extensive, training operator learning models can also be challenging. 

For S2 and S3, the core issue involves geometric encoding issues of collocation points and the computation domain. There are already several existing solutions. The approaches broadly fall into two categories: 
\begin{enumerate}
    \item Mathematics perspective: Transforming the original problem into an equivalent formulation \cite{WOS:000612233300006}. Such as, Minglang Yin proposed Diffeomorphic Mapping Operator Learning (DIMON) \cite{WOS:001372682600001}, which incorporates diffeomorphisms with DeepONet to solve geometry-dependent problems. Meng B. et al. proposed three operator learning models, based on boundary integral equations, to solve 2D PDEs where the boundary is smooth and can be expressed by boundary integral equations \cite{WOS:001438180600001}. Zongyi Li proposed GeoFNO, which transforms computational domains into regular geometries via neural networks before learning with FNO\cite{JMLR:v24:23-0064}. Shanshan Xiao solved 2D Poisson-related problems by integrating transformations with MIONET \cite{Xiao2024LearningSO}.
    \item Computer perspective: Leveraging techniques like point clouds \cite{KASHEFI2022111510}, signed distance functions \cite{WOS:001256289700001}, advanced neural networks and training technique \cite{10.5555/3666122.3666393, WOS:001049372700001, tran2023factorized, UnisolverHZhou, 10.5555/3666122.3669215}, and transfer learning \cite{WOS:001147270900002} and so on. For instance, Chenyu Zeng proposed the Point Cloud Neural Operator (PCNO), combining point clouds with neural operators \cite{WOS:001498427300001}. Zongyi Li and associates proposed the geometry-informed neural operator (GINO), based on SDF, point clouds, and neural operators using graph and fourier architectures, to solve 3D flow problems around vehicles \cite{NEURIPS2023_70518ea4}. Zhongkai Hao proposed GNOT \cite{10.5555/3618408.3618917}, a transformer-based and mixture-of-experts-based (MoE) architecture, to address multiple related PDEs. Ning Liu embeds computational domains into rectangular regions, introducing the Domain Agnostic Fourier Neural Operator (DAFNO) to learn surrogates for irregular geometries and evolving domains \cite{liu2023domain}. 
\end{enumerate}
For the former, establishing an equivalent formulation—whether through computational domain transformations or boundary integral equations—poses significant challenges. For the latter, the cost involves substantial computational overhead, particularly pronounced for point cloud techniques when scaling to datasets with massive samples and grid points \cite{Transolver++Luo}.

We introduce the three-type DD-DeepONet, two iteration frameworks, and an iteration-free framework, which integrates DDM with DeepONet, addressing the requirements of three application scenarios to a certain extent. Combined with stretching transformations and translation, it can be used to solve geometry-dependent problems. This method combines the advantages of DDM and operator learning by employing DeepONet as a standalone solver. This framework also features two stages: off-line (training) and on-line (inference) stages. Data preparation and computational strategies differ across stages depending on the iteration scheme selected.

It is important to note that this is a conceptual approach, not confined to specific iterative formats or solvers. The iterative scheme and DeepONet structure utilized in this paper can be substituted, and subdomains may be solved using any numerical method. For instance, Minglang Yin has integrated finite elements with DeepONet to solve multi-scale mechanical problems \cite{WOS:000897844700005}, while the combination of traditional numerical methods and DDM has been widely adopted. We chose the DeepONet and MIONET structure here because we need to handle multiple function spaces.

We solve PDEs in three application scenarios. S1: 3D Laplace equation with mixed boundaries, S2: 3D Laplace for resistance, 2D steady Navier-Stokes (N-S) equations for pipe flow, S3: 2D steady drift-diffusion equations for charge carriers, 3D multimedium  Poisson equation.

The remainder of this paper is organized as follows. In Section 2, we provide a concise mathematical description of the operators to be learned in this paper. In Section 3, we briefly review DDM and the relevant knowledge of DeepONet, and introduce DD-DeepONet. Section 4 presents several numerical examples. Finally, in Section 5, we conclude and discuss the paper and provide an outlook for future research.

\section{Description of the problem}
Here, we provide a mathematical general framework description of the mappings required for three application scenarios. 

Let $\Omega = \cup_{i} D_{i} \subset \mathbb{R}^{d}$ be a simply connected, closed, and bounded domain, composed of the union of various regions $D_{i}$. Each $D_{i}$ represents a subdomain geometric structure, such as a cuboid or rectangle, in this paper. Let $\mathscr{L}_{E}$ denote a partial differential operator, and $\mathscr{L}_{B}$  BCs operator, we consider the following problem:
\begin{equation}
\begin{aligned}
    (\mathscr{L}_{E} (\alpha, u))(\bm{x}) &= p(\bm{x}), \quad \bm{x} \in \Omega, \\
    (\mathscr{L}_{B} (\beta, u))(\bm{x}) &= q(\bm{x}), \quad \bm{x} \in \partial \Omega,
\end{aligned}
\label{eq_operator}
\end{equation}
for $\alpha, \beta \in \mathbb{R}, p(\bm{x}) \in \mathcal{P}, q(\bm{x}) \in \mathcal{Q}$, where $\alpha, \beta$ represent the coefficients of equation, $\mathcal{P}, \mathcal{Q}$ denote Banach spaces. $u$ represents the PDE's solution that needs to be solved. 

We assume that $\mathscr{L}_{E}$ and $\mathscr{L}_{B}$ are such that, for any quintuple $(\Omega, \alpha, \beta, p, q)$, the PDE (\ref{eq_operator}) has unique solution $u(\bm{x}) \in \Phi$ where $\Phi$ denotes a Banach space defined on $\Omega$. Consequently, the solution to the equation can be expressed as 
\begin{equation}
    \phi(\bm{x}) = \mathcal{G}[\Omega, \alpha, \beta, p, q](\bm{x}), \quad  \bm{x} \in \Omega. 
\label{eq_operator_sol}
\end{equation}

We can define the mapping from the product space to the solution space
\begin{equation}
    \mathcal{G}: \mathbb{R}^{d} \times \mathbb{R} \times \mathbb{R} \times \mathcal{P} \times \mathcal{Q} \rightarrow \Phi. 
\label{eq_operator_map}
\end{equation}
The operator $\mathcal{G}$, which has been abstracted, represents the operator we ultimately aim to approximate. It maps between infinite-dimensional function spaces. 

Depending on the specifics of the problem, the number of product spaces might need adjustment, and they are not necessarily composed of these five Banach spaces. 

For Scenario 1, the mapping relationship that needs to be learned can be denoted as $\mathcal{G}: \mathcal{R} \times \mathcal{R} \times \mathcal{P} \times \mathcal{Q} \rightarrow \Phi$. 

For Scenario 2, the required mapping to be learned is represented as $\mathcal{G}: \mathbb{R}^{d} \rightarrow \Phi$. 

For Scenario 3, i.e. more general situation, is mapping (\ref{eq_operator_map}).

These mappings are crucial for understanding how various aspects of a problem's setup—such as geometry, BCs, and physical parameters—affect the solution space. This understanding, in turn, allows for the application of suitable computational strategies to address a wide range of engineering and scientific challenges. 

\section{Method}
In this section, we will briefly introduce DDM and the basics of DeepONet. Then, we introduce our proposed framework, DD-DeepONet. 

\subsection{Domain decomposition method}
DDM encompasses a suite of techniques based on the divide-and-conquer principle. To date, numerous DDMs have been developed, such as the Schwarz alternating method, finite element tearing and interconnecting (FETI), multigrid methods, and preconditioning methods induced by Schur complement, among others \cite{Saad_DDM, AT_OBW_DDM, WOS:000169338900003}. 

\begin{figure}[htbp]
\centering
\includegraphics[scale=0.25]{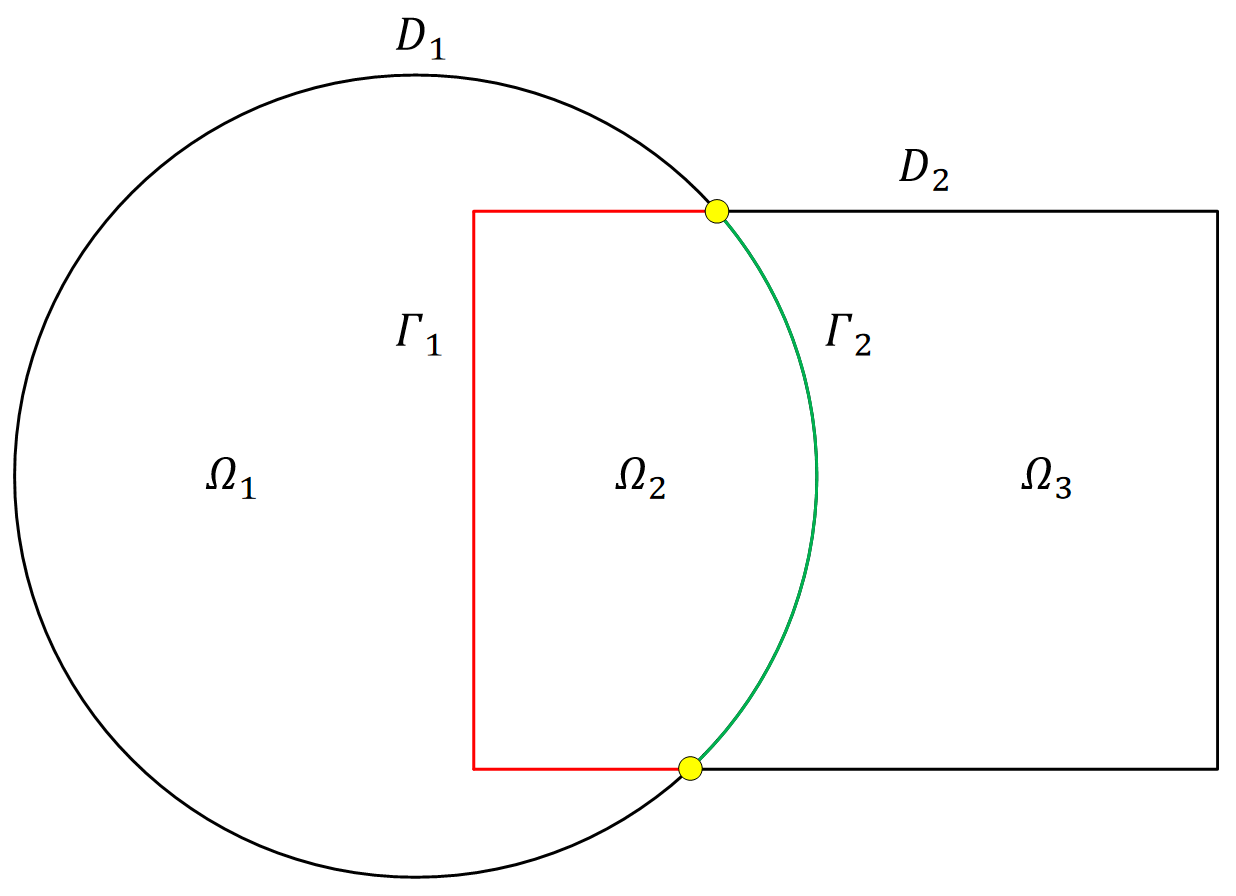}
\caption{Overlap partition for the Schwarz alternating method with two subdomains.}
\label{fig3.1.1}
\end{figure}

Let's review the earliest known Schwarz alternating method \cite{WOS:000086935000018, AT_OBW_DDM}. The description of the two-subdomain Schwarz alternating method is as follows: Consider PDEs defined on a bounded Lipschitz domain $\Omega = \Omega_{1} \cup \Omega_{2} \cup \Omega_{3}$,
\begin{equation}
\begin{aligned}
    - \Delta u & = f, \quad in \, \Omega, \\
    u & = 0, \quad on \, \partial \Omega. 
\end{aligned}
\label{eq3. 1. 1}
\end{equation}
The subdomians $D_{1} = \Omega_{1} \cup \Omega_{2}$ and $D_{2} = \Omega_{2} \cup \Omega_{3}$ are circular and rectangular areas, respectively. The geometric description of these regions is illustrated in Figure \ref{fig3.1.1}. 

\begin{algorithm}
\caption{The Schwarz alternating method. }\label{Algorithm_1}
\begin{algorithmic}
\State \textbf{Initialize:} Set $u^{0}=0$ on $\Gamma_{2}$. 
\State \textbf{Main Loop:}
\For{$n = 0 : N$}
    \State \textbf{Circular domian:}
            \begin{itemize}[label=\textbullet, font=\color{black}\large]
              \item Receive interface information $u^{n}$ on $\Gamma_{2}$ from cuboid domain, or from initial guass $u^{0}$. 
              \item Solve for $u^{n + 1/2}$ in $D_{1}$ from following equation:
                    \begin{equation}
                    \begin{aligned}
                        - \Delta u^{n+1/2} & = f, \quad in \, D_{1}, \\
                        u^{n+1/2} & = 0, \quad on \, \partial  D_{1} \textbackslash \Gamma_{2}, \\
                        u^{n+1/2} & = u^{n}, \quad on \, \Gamma_{2},
                    \end{aligned}
                    \label{eq3.1.2}
                    \end{equation}
              \item Compute and collect interface information $u^{n + 1/2}$ on $\Gamma_{1}$ and pass it to cuboid domian. 
            \end{itemize}
    \State \textbf{Rectangular domian:}
            \begin{itemize}[label=\textbullet, font=\color{black}\large]
              \item Receive interface information $u^{n + 1/2}$ on $\Gamma_{1}$ from circular domain. 
              \item Solve for $u^{n + 1}$ in $D_{2}$ from following equation:
                    \begin{equation}
                    \begin{aligned}
                        - \Delta u^{n+1} & = f, \quad in \, D_{2}, \\
                        u^{n+1} & = 0, \quad on \, \partial  D_{2} \textbackslash \Gamma_{1}, \\
                        u^{n+1} & = u^{n + 1/2}, \quad on \, \Gamma_{1},
                    \end{aligned}
                     \label{eq3.1.3}
                    \end{equation}
              \item Compute and collect interface information $u^{n + 1}$ on $\Gamma_{2}$ and pass it to circular domian. 
            \end{itemize}
    
    \State \textbf{If converged, stop;}
\EndFor
\end{algorithmic}
\end{algorithm}

Algorithm \ref{Algorithm_1} summarizes the iterative procedure of the original Schwarz alternating method. At the start of the iteration, an initial guess of $u^{0}=0$ is given on $\Gamma_{2}$, followed by solving the Laplace boundary value problem (\ref{eq3.1.2}) to obtain the solution $u^{n + 1/2}$ in the circular computational region $D_{1}$. Here, the original BCs of problem (\ref{eq3. 1. 1}) are maintained on the boundary $\partial D_{1} \setminus \Gamma_{2}$, and the initial guess is used on the newly constructed fictitious boundary $\Gamma_{2}$. Next, interface information $u^{n + 1/2}$ on $\Gamma_{1}$ is calculated and collected within the circular region $D_{1}$, and this information is passed to the rectangular region. In the rectangular area $D_{2}$, the Laplace boundary value problem (\ref{eq3.1.3}) is solved to obtain its solution $u^{n+1}$, where the original BC of problem (\ref{eq3. 1. 1}) are maintained on the boundary $\partial D_{2} \setminus \Gamma_{1}$, and the values transferred from the circular region $D_{1}$ are used on the new fictitious boundary $\Gamma_{1}$. Interface information $u^{n+1}$ is then calculated and collected on $\Gamma_{2}$ within the rectangular computational area $D_{2}$, and this information is passed back to the circular region. The iteration is incremented from $n$ to $n+1$, and the steps are repeated cyclically until the iteration meets the specified error criteria and is subsequently terminated. Additionally, to accelerate the convergence of the solution process, a relaxation factor $\theta$ can be introduced when updating the interface information, $\bar{u}^{n+1} = (1-\theta) \cdot u^{n} + \theta \cdot u^{n + 1}$. 

\subsection{DeepONet}

\begin{figure}[htbp]
\centering
\includegraphics[scale=0.45]{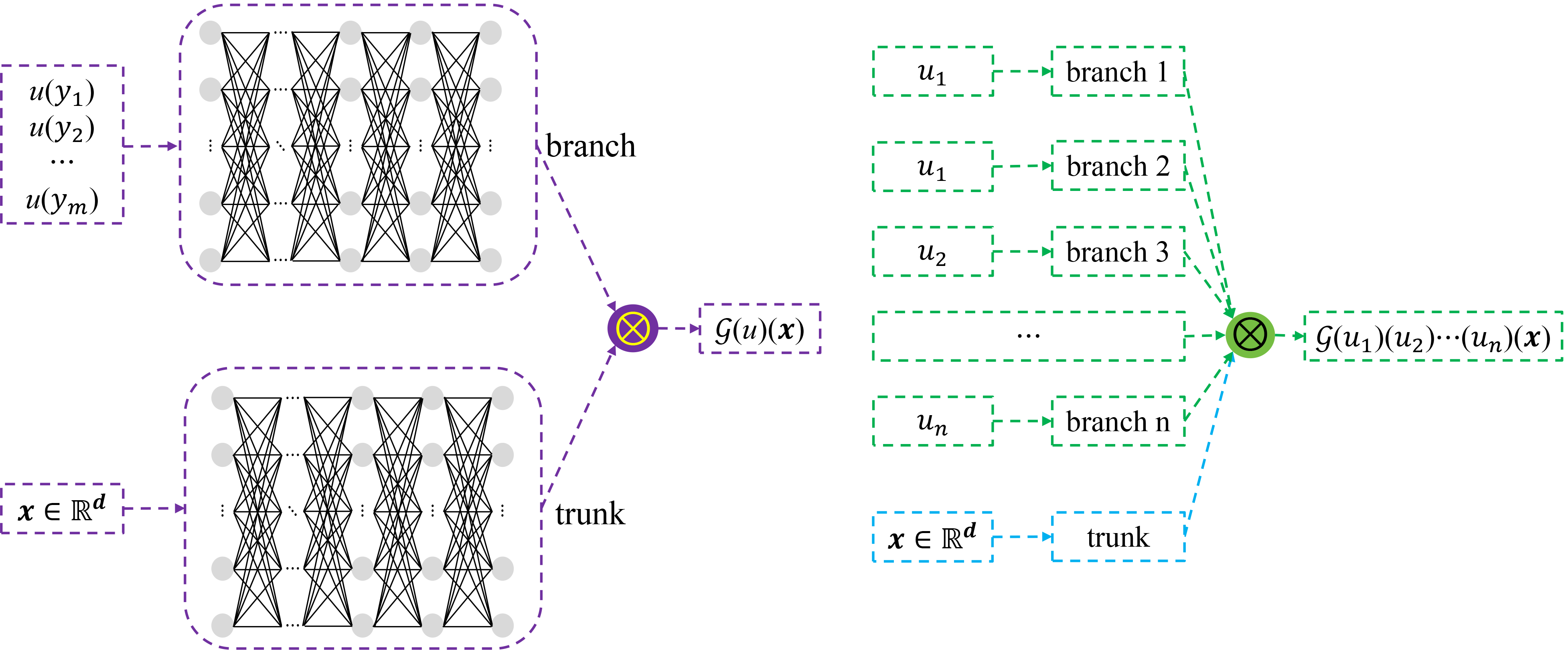}
\caption{(left) Vanilla DeepONet. (right) MIONET or EDeepONet.}
\label{fig3.2.1}
\end{figure}

Following the universal approximation theorem for functions, Chen and others introduced a universal approximation theorem for operators in 1995 \cite{392253}. Building on this foundation, Lulu and colleagues further developed and proposed the DeepONet structure \cite{WOS:000641834300001}, see Figure \ref{fig3.2.1} (left), which generalizes the universal approximation theorem for operators. This network architecture can be used to learn both explicit and implicit operators, including nonlinear ones, and has demonstrated remarkable generalization capabilities. Pengzhan Jin and Lesley Tan proposed the multiple-input operator network (MIONET) \cite{WOS:000885480600003} and enhanced DeepONet \cite{tan2022enhancedDeepONetmodelingpartial} (EDeepONet), respectively, both of which extend to multiple input Banach product space, as shown in Figure \ref{fig3.2.1} (right). In the following parts, we provide a brief overview of the original DeepONet. 

Suppose that $\Omega \subset \mathbb{R}^{d}$ is a bounded open set, which serves as the domain for PDE problems. In this paper, we focus solely on issues related to PDEs. The original DeepONet learns the mapping between two Banach spaces, which take values in $\mathbb{R}^{m}$ and $\mathbb{R}^{d}$, respectively. We denote these two Banach spaces as $\mathcal{U} = \mathcal{U} (\Omega_{branch};\mathbb{R}^{m})$ and $\mathcal{S} = \mathcal{S}(\Omega;\mathbb{R}^{d})$. The mapping approximated by DeepONet is $\mathcal{G}: \mathcal{U} \rightarrow \mathcal{S}$, meaning that given $u \in \mathcal{U}$, the corresponding solution to the PDEs is $\mathcal{G}(u) \in \mathcal{S}$. Here, the Banach space $\mathcal{U}$ may represent the space of BC or parameters in the equations, while $\mathcal{S}$ is the space of solutions to the equations. For any $\bm{x} \in \mathbb{R}^{d}$, $\mathcal{G}(u)(\bm{x}) \in \mathbb{R}$ represents the evaluation of the function $\mathcal{G}(u)$ at $\bm{x}$:
\begin{equation}
    \mathcal{G}(u)(\bm{x}) \approx branch \cdot trunk.
\end{equation}
In this work, we adopt FCNN as the sub-architecture.

\subsection{DD-DeepONet}

\begin{figure}[htbp]
\centering
\includegraphics[scale=0.5]{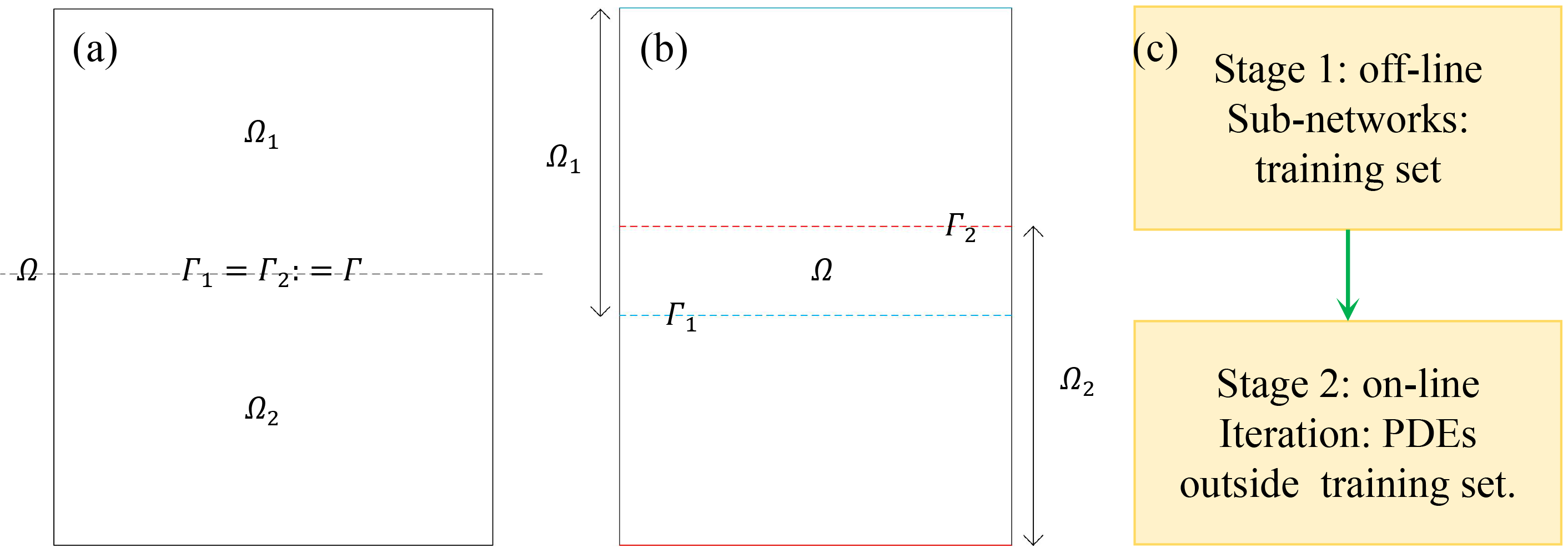}
\caption{DD-DeepONet Computational Domain and Workflow. (a) Non-overlap. (b) Overlap. (c) Workflow.}
\label{fig3.3.1}
\end{figure}

In this section, we introduce DD-DeepONet for PDEs. DD-DeepONet combines DDM with DeepONet, using DeepONet as a solver for subdomains. We will use a 2D two-subdomain DDM as an example to elucidate our two computational frameworks, as shown in Figure \ref{fig3.3.1}(a) and (b).

We note that this framework can not only be extended to multiple subregions and other dimensions, but its iterative framework can also be replaced to yield alternative approaches. Additionally, subdomain networks may be substituted with traditional numerical methods or other complex networks, demonstrating great flexibility.

This framework consists of two stages, Figure \ref{fig3.3.1}(c): Stage 1 (off-line): Determine the iteration scheme, collect subdomain training data, and train/save subdomain networks. Subdomain datasets can be obtained by: Processing global data within subdomains locally; Using subdomain prior knowledge to approximate or fit the function space of BCs, then solving subdomain PDEs to generate data. Stage 2 (on-line): For PDEs outside the training set, solve them via the coupling framework (detailed next), assuming completion of the first stage with saved subdomain networks.

\subsubsection{The coupling framework 1}

\begin{figure}[htbp]
\centering
\includegraphics[scale=0.55]{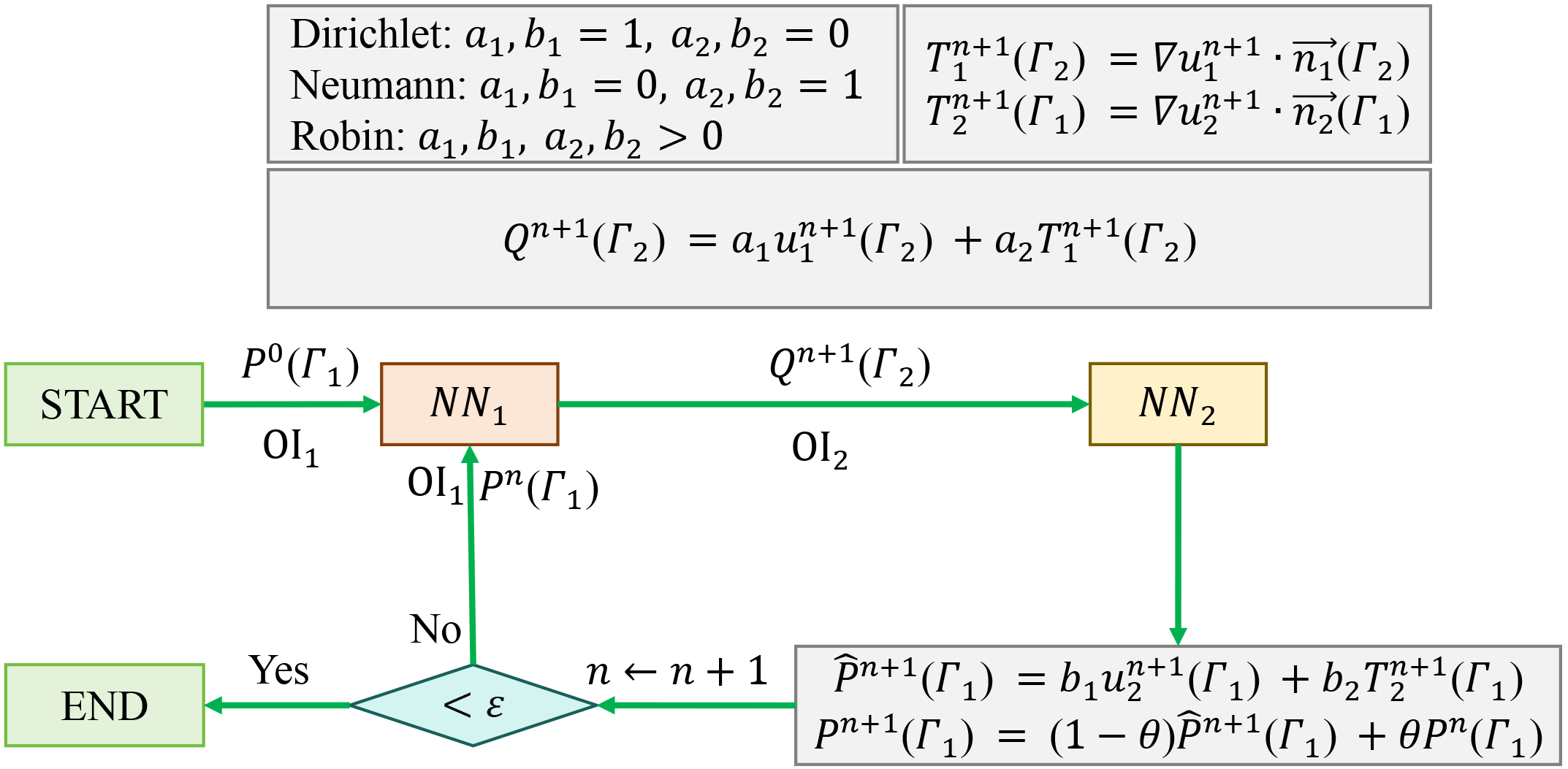}
\caption{The coupling framework 1. $T_{i}^{n+1}$ denotes the outward normal derivative, $\theta$ is a relaxation factor, and $OI_{i}$ refers to other inputs of neural networks. The parameters $a_{i}$ and $b_{i}$ can take different values, allowing $P^{n+1}$ and $Q^{n+1}$ to impose various types of BC, so it can form different iterative schemes.}
\label{fig3.3.2}
\end{figure}

\begin{algorithm}
\caption{The coupling framework1. }\label{Algorithm_2}
\begin{algorithmic}
\State \textbf{Initialize:} Select a proper $P^{0}$ on $\Gamma_{1}$, prepare $OI_{i}, i = 1, 2$. 
\State \textbf{Main Loop:}
\For{$n = 0 : N$}
    \State \textbf{$NN_{1}$ in $\Omega_{1}$:}
            \begin{itemize}[label=\textbullet, font=\color{black}\large]
              \item Receive interface information $P^{n}$ on $\Gamma_{1}$ from $NN_{2}$, or from initial guass $P^{0}$. 
              \item Input $P^{n}$ and $OI_{1}$ into $NN_{1}$ to compute the solution $u_{1}^{n+1}$ over the domain $\Omega_{1}$. 
              \item Calculate $Q^{n+1} = a_{1} u_{1}^{n+1} + a_{2} T_{1}^{n+1}$ on $\Gamma_{2}$ and pass it to $NN_{2}$. 
            \end{itemize}
    \State \textbf{$NN_{2}$ in $\Omega_{2}$:}
            \begin{itemize}[label=\textbullet, font=\color{black}\large]
              \item Receive interface information $Q^{n+1}$ on $\Gamma_{2}$ from $NN_{1}$. 
              \item Input $Q^{n+1}$ and $OI_{2}$ into $NN_{2}$ to compute the solution $u_{2}^{n+1}$ over the domain $\Omega_{2}$. 
              \item Calculate $\hat{P}^{n+1} = b_{1} u_{2}^{n+1} + b_{2} T_{2}^{n+1}$ and $P^{n+1}=(1-\theta)\hat{P}^{n+1} + \theta P^{n}$ on $\Gamma_{1}$ and pass it to $NN_{1}$. 
            \end{itemize}
    
    \State \textbf{If converged, stop;}
\EndFor
\end{algorithmic}
\end{algorithm}

Figure \ref{fig3.3.2} and Algorithm \ref{Algorithm_2} describe the algorithmic process of coupling framework 1. Applicable to both overlapping and non-overlapping cases. At the start of iteration, an initial guess $P^{0}$ is given on $\Gamma_{1}$, and other inputs for $NN_{1}$ and $NN_{2}$ are collected. Once the iteration begins ($n \geq 0$), $P^{0}$ and $OI_{1}$ are fed into $NN_{1}$ to compute the neural network's output $u_{1}^{n+1}$ on $\Omega_{1}$. Subsequently, calculate $T_{1}^{n+1}$ on interface $\Gamma_{2}$, if necessary. Several methods could be employed here, such as finite difference or neural network. Next, we compute the inputs $Q^{n+1} = a_1 u_{1}^{n+1} + a_{2} T_{1}^{n+1}$ for $NN_{2}$ on $\Gamma_{2}$, and feed them along with $OI_{2}$ into $NN_{2}$ to compute the output $u_{2}^{n+1}$ on $\Omega_{2}$. Using similar methods, $T_{2}^{n+1}$ on $\Gamma_{1}$ can be computed, and then $\hat{P}^{n+1} = b_{1} u_{2}^{n+1} + b_{2} T_{2}^{n+1}$ can be determined. Then update $P^{n+1}$ by $P^{n+1} = (1-\theta)\hat{P}^{n+1} + \theta P^{n}$, where the relaxation factor $\theta \in [0,1]$, is either fixed at a value or updated dynamically. If the termination criteria are met, the DDM is considered converged, and the iteration stops. If not, the iteration continues by setting $P^{n+1}$ as $P^{n}$ and re-entering the loop for further calculation. 

There are many termination criteria that can be chosen, such as specifying a number of iterations $iter < C$, or defining the difference between two successive iterations,
\begin{equation*}
    \| u_{1}^{n+1} - u_{1}^{n} \| + \| u_{2}^{n+1} - u_{2}^{n} \| < \epsilon. 
\end{equation*}

\subsubsection{The coupling framework 2}

\begin{figure}[htbp]
\centering
\includegraphics[scale=0.55]{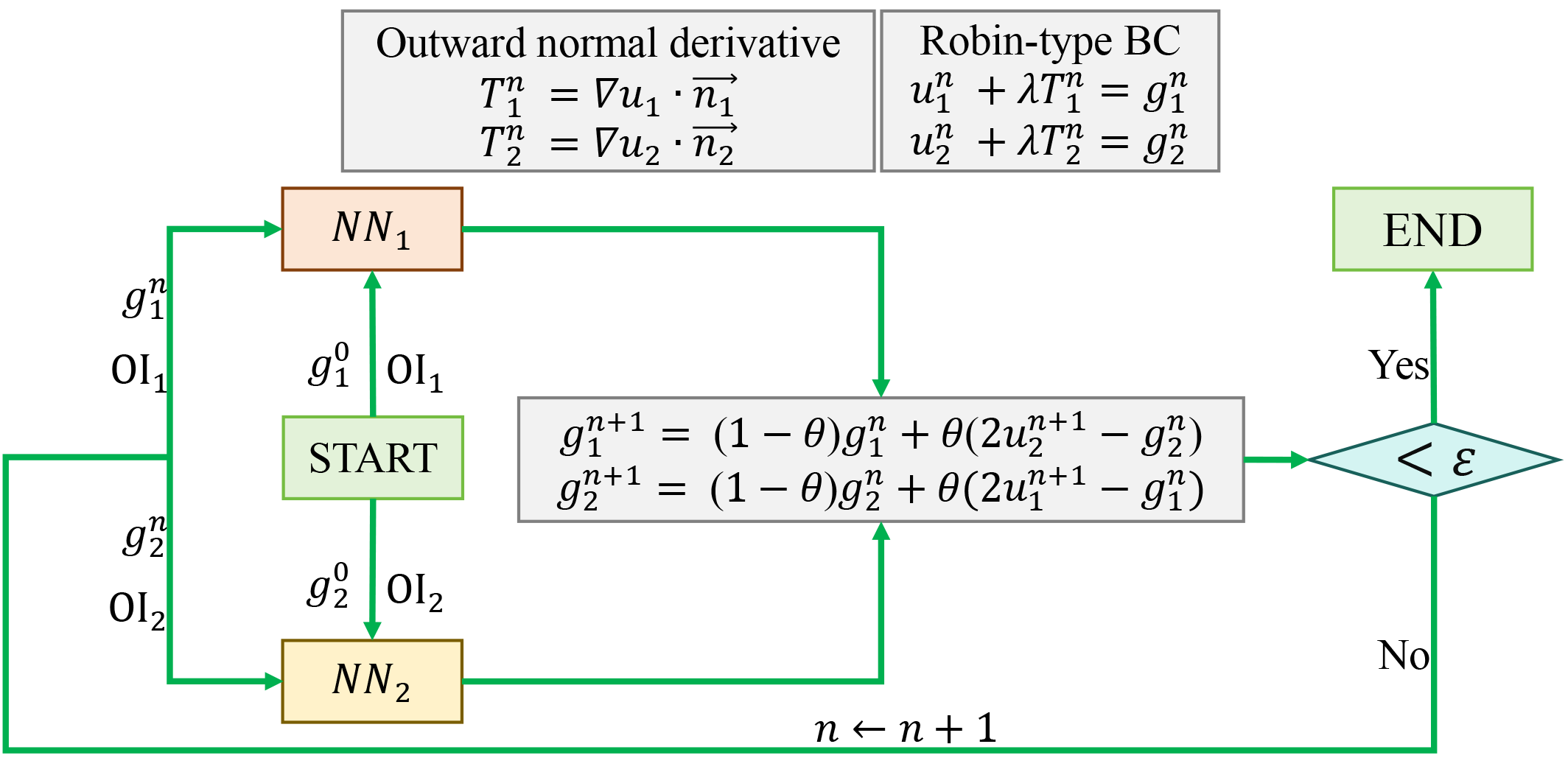}
\caption{The coupling framework 2. $g_{i}^{n}$ represents the Robin-type interface transmission condition. Other notations are similar to those in Figure \ref{fig3.3.2}.}
\label{fig3.3.3}
\end{figure}

\begin{algorithm}
\caption{The coupling framework2. }\label{Algorithm_3}
\begin{algorithmic}
\State \textbf{Initialize:} Select proper $g_{i}^{0}$ on $\Gamma$, prepare $OI_{i}, i = 1, 2$. 
\State \textbf{Main Loop:}
\For{$n = 0 : N$}
    \State \textbf{$NN_{1}$ in $\Omega_{1}$:}
            \begin{itemize}[label=\textbullet, font=\color{black}\large]
              \item Receive interface information $g_{1}^{n}$ on $\Gamma$, or from initial guass $g_{1}^{0}$. 
              \item Input $g_{1}^{n}$ and $OI_{1}$ into $NN_{1}$ to compute the solution $u_{1}^{n+1}$ over the domain $\Omega_{1}$ and $u_{1}^{n+1}$ on $\Gamma$. 
            \end{itemize}
    \State \textbf{$NN_{2}$ in $\Omega_{2}$:}
            \begin{itemize}[label=\textbullet, font=\color{black}\large]
              \item Receive interface information $g_{2}^{0}$ on $\Gamma$, or from initial guass $g_{2}^{0}$. 
              \item Input $g_{2}^{n}$ and $OI_{2}$ into $NN_{2}$ to compute the solution $u_{2}^{n+1}$ over the domain $\Omega_{2}$ and $u_{2}^{n+1}$ on $\Gamma_{1}$. 
            \end{itemize}
    \State \textbf{Updata Robin-type BCs:}
            \begin{itemize}[label=\textbullet, font=\color{black}\large]
              \item Calculate $g_{1}^{n+1}=(1 - \theta) g_{1}^{n} + \theta (2 u_{2}^{n+1} - g_{2}^{n})$, pass it to $NN_{1}$. 
              \item Calculate $g_{2}^{n+1}=(1 - \theta) g_{2}^{n} + \theta (2 u_{1}^{n+1} - g_{1}^{n})$, pass it to $NN_{2}$. 
            \end{itemize}
    \State \textbf{If converged, stop;}
\EndFor
\end{algorithmic}
\end{algorithm}

Figure \ref{fig3.3.3} and Algorithm \ref{Algorithm_3} describe the computational process of the coupling framework 2, which is used for the non-overlapping case. The domain is shown in Figure \ref{fig3.3.1} (a). This approach combines Lions' DDM with DeepONet. 

Lions' DDM is a method for non-overlapping DDM that uses Robin BCs at the interfaces between subdomains to exchange data \cite{WOS:000168784700021, INSPEC:5724898}. Notably, its iterative scheme does not require the calculation of normal derivatives, which can significantly reduce the computational load and minimize errors caused by calculating derivatives, thereby greatly enhancing convergence. Furthermore, this method allows for significant parallelization.

At the start of the iteration, initial guesses for the two Robin-type boundary transmission conditions defined on $\Gamma$ are given as $g_{i}^{0}$, which, along with the other neural network inputs $OI_{i}$, are fed into the trained network models $NN_{1}$ and $NN_{2}$ to compute the solutions $u_{i}^{n+1}$ in their respective subdomains. Then update the two Robin-type BCs $g_{i}^{n+1}$. Subsequently, check if the termination criteria are met; if so, the iteration stops; otherwise, it sets $n+1$ as $n$ and continues the update steps until convergence. 

\subsubsection{Iteration-free}

Traditional DDM and the previous two coupling frameworks not only divide the computational domain into multiple subdomains but also require the exchange of information at interfaces to drive iteration, demanding strict continuity at these interfaces. For the DeepONet structure, the trunk net inputs are the coordinate points within the computational domain. If we do not take into account the continuity of the contact part, we can fully embrace DDM by training multiple sub-networks. Each network inputs coordinate points from its respective computational domain, without exchanging information and iteration. But each branch net needs to input the global information. 

\section{Numerical examples}

In this section, we present several numerical examples calculated using our method: DD-DeepONet. All simulation data were conducted with Intel(R) Xeon(R) Platinum 8380 CPU. All neural network training and inference were conducted using a single Nvidia A100 GPU card. 

The setup of the dataset and neural network is shown in \ref{Appendix_A} (Table \ref{table_NN_Setup_FEx}, \ref{table_NN_Setup_FEx_resistance}, \ref{table_NN_Setup_FEx_pipeflow}, \ref{table_NN_Setup_FEx_DD}, \ref{table_NN_Setup_possion}). The neural networks are all initialized using the Xavier initialization, with ReLU as the activation function and Adam as the optimizer. The batch size is 64. The learning rate is optimized using an exponential decay strategy, initialized at 0.0002. 

In this paper, we utilize the mean $L_{2}$ relative error (ML2RE) as the loss function:
\begin{equation}
    ML2RE = \frac{1}{N} \sum^{N}_{i=1} \frac{\| \bm{y}_{i}^{pred} - \bm{y}_{i}^{true} \|_{2}}{\| \bm{y}_{i}^{true} \|_{2}} = \frac{1}{N} \sum_{i=1}^{N} \frac{\sqrt{\sum_{j=1}^{M}(y_{ij}^{pred} - y_{ij}^{true})^{2}}}{\sqrt{\sum_{j=1}^{M} (y_{ij}^{true})^{2}}},
\label{eq3.2.3}
\end{equation}
$\bm{y}_{i}^{pred}$ represents the predicted values, while $\bm{y}_{i}^{true}$ denotes the true values for the $i$-th sample respectively. $N$ denotes the number of samples, and $M$ represents the number of collocation points per sample. 

\subsection{S1 examples}

\begin{figure}[htbp]
\centering
\includegraphics[scale=0.45]{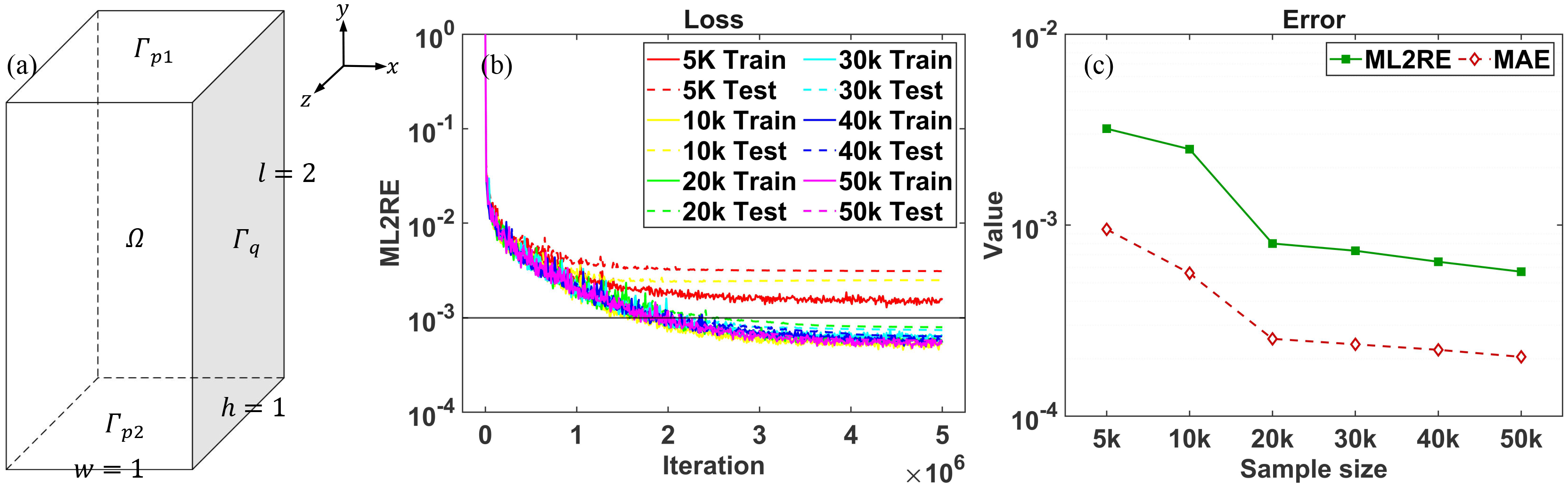}
\caption{Computational domain and sample size impact. (a) Cuboid computational domain of the Laplace equation. (b) Loss curve. (c) Variation of generalization error with sample size.}
\label{fig4.1.1}
\end{figure}

Consider the problem of solving a 3D Laplace equation, where the computational domain is defined within a cuboid, see Figure \ref{fig4.1.1} (a). The equation satisfied is as follows:
\begin{equation}
   \begin{aligned}
    - \Delta u & = 0, \quad in \, \Omega, \\
    u & = f, \quad on \, \Gamma_{p1} \cup \Gamma_{p2}, \\
    \frac{\partial u}{\partial \vec{n}} &= 0, \quad on \, \Gamma_{q}. 
\end{aligned}
\label{eq4.1.1}
\end{equation}
$u$ is the solution to the equation. $\Gamma_{pi}$, $i=1,2$, represent non-homogeneous Dirichlet BCs, and $\Gamma_{q}$ represents homogeneous Neumann BCs imposed on the four lateral faces of the cuboid. 

The generation of Dirichlet BCs $f$ is achieved using Gaussian Process (GP \cite{WOS:000641834300001, lu2021deepxde}), with the correlation length set to $l=0. 5$. We use our custom-written $C++$ finite element code to generate data. All samples were discretized using second-order 27-node hexahedral Lagrange elements. Each sample has 99937 grid points. All numerical results were saved in VTK format files. We got 51,000 samples.

The presented computational results are derived from a generalization set of 1,000 samples, which was explicitly excluded from both the training and test datasets. This dedicated validation subset was specifically designed to evaluate how DDM affects model generalization capabilities.

We use ML2RE and mean absolute error (MAE) to characterize the error:
\begin{equation}
    MAE = \frac{1}{N} \sum_{i=1}^{N} |y_{i}^{true} - y_{i}^{pred}| = \frac{1}{NM} \sum_{i=1}^{N} \sum_{j=1}^{M} |y_{ij}^{true} - y_{ij}^{pred}|, 
\label{eq4.1.2}
\end{equation}
$N$ denotes the number of samples, and $M$ represents the number of collocation points per sample.

We first tested the impact of dataset size on model performance. We selected subsets of 5,000, 10,000, 20,000, 30,000, 40,000, and 50,000 samples from a dataset of 51,000 samples, each forming a new dataset.  Figure \ref{fig4.1.1}(b) and (c) show that the error reduction rate decreases significantly when the sample size reaches 20,000.

This section considers three distinct and representative DD-DeepONet schemes.

\subsubsection{Iteration-free}

\begin{figure}[htbp]
\centering
\includegraphics[scale=0.65]{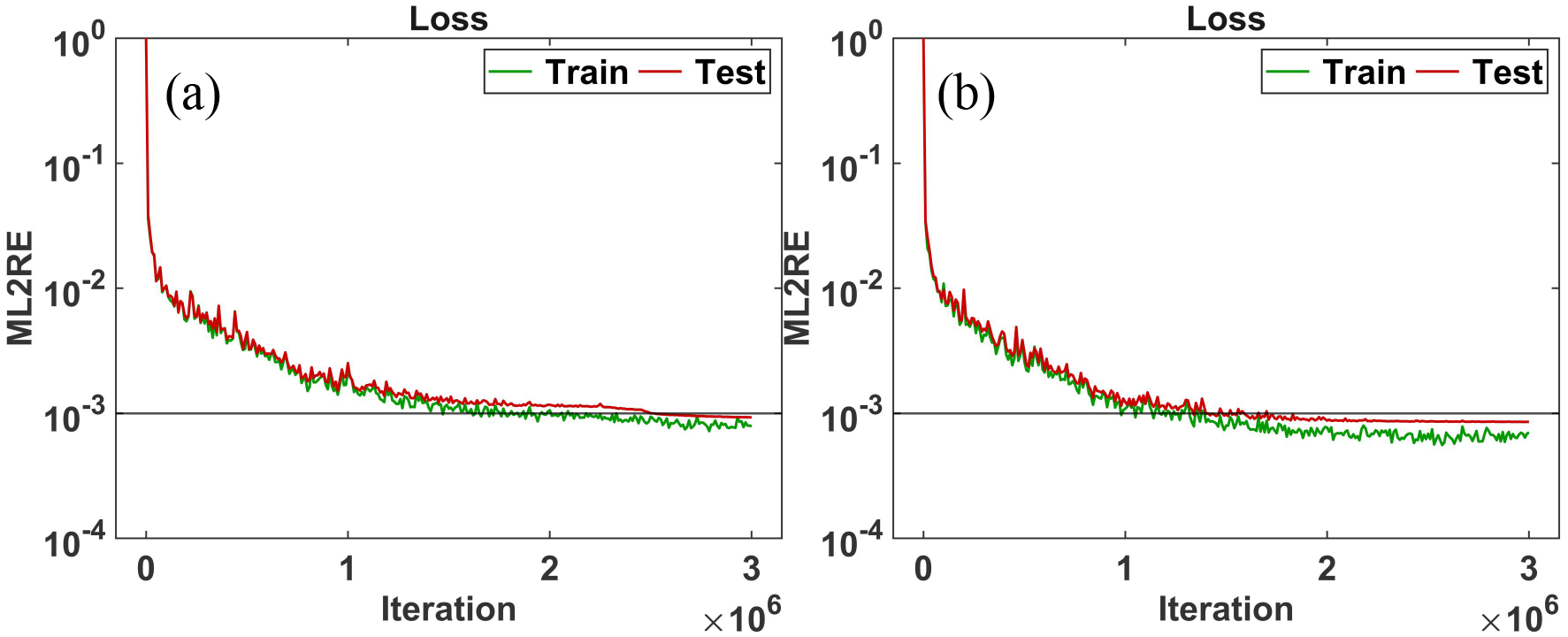}
\caption{Comparison of non-DDM and DDM results. (a) Loss curves without DDM. (b) Loss curves using DDM. }
\label{fig_FEx_case1}
\end{figure}

\begin{table}[pos=htbp]
\centering
\caption{Comparison between DDM and non-DDM results.} 
\begin{tabularx}{\textwidth}{CCCCC} 
\toprule
 & ML2RE & MAE & TIME & MEMORY \\
\hline
     Non-DDM & 9.2951e-4 & 2.4089e-4 & 344,106$s$ & 15,217$M$  \\
     DDM(half)& 8.7841e-4 & 1.8019e-4 & 206,240$s$ & 8,575$M$  \\
\bottomrule
\end{tabularx}
\label{table_Iteration-free}
\end{table}

We train neural networks on two domains for comparison: the global domain ($0 \leq y \leq 2$) without DDM, and a subdomain ($0 \leq y \leq 1$) using half collocation points. The only difference is the number of collocation points input to the trunk network, while all other settings remain identical. Figure \ref{fig_FEx_case1} and Table \ref{table_Iteration-free} show that the iteration-free DD-DeepONet alleviates training difficulty, shortens training time, eases GPU memory usage, and slightly improves accuracy.

Although global domain learning requires training only a single neural network, it is well-known that as the computational scale increases, the accuracy of a fixed-size network tends to decrease due to overfitting. While enlarging the network size can theoretically enhance accuracy, this approach does not always hold; excessively large networks lead to severe training difficulties and lower accuracy. DDM effectively addresses these issues \cite{WOS:001147270900002}. Although DD-DeepONet requires training multiple subnetworks for global computations, it delivers significant overall advantages.

\subsubsection{D-R (non-overlap) and D-D (overlap)}

In this subsection, we implement two iteration schemes from framework 1—non-overlapping and overlapping—for equation (\ref{eq4.1.1}): D-R and D-D type DD-DeepONet. We consider two data acquisition methods for subdomains:
 \begin{enumerate}
    \item Perform interpolation on the global domain dataset to collect subdomain training data;
    \item Empirically assume the function space of interface data. Generate substantial interface conditions, supplement missing BCs for subdomain problems to ensure well-defined. Then, solving PDEs that are defined on subdomains, collect subdomain training data.
\end{enumerate}
For both D-R and D-D schemes, we analyze and compare the two methods using interpolation and GP($l=0.5$), respectively. 

Based on preliminary tests where errors stabilized beyond 20,000 samples, we set all dataset sizes to 30,000 to balance computational costs and ensure a fair comparison. For DD-DeepONet iterations, we adopt: convergence criterion is $MAE < 1e-4$; relaxation factor is $\theta = 0.5$. The initial iteration values at the interface are all set to 0. For both iteration schemes, the computational domain is partitioned into two subdomains along the y-axis. Each D-R subdomain spans a length of 1, while each D-D subdomain spans 1.25.

The D-R scheme implements Robin BC: $Q = u + \nabla u \cdot \vec{n}$, with $\nabla u$ computed via second-order finite differences on equidistant mesh nodes:
\begin{equation}
    \dot{u}(x) = \frac{-4 u(x-h) + u(x-2h) + 3u(x)}{2h}.
\label{eq4.1.4}
\end{equation}

In the D-D scheme, we observe that the two subdomains ($0 \leq y \leq 1.25$ and $0.75 \leq y \leq 2$) are identical under a 0.75-unit y-axis translation. Thus, we merge datasets that are generated by GP and train a single neural network, reducing the required subdomain networks and improving accuracy. 

\begin{figure}[htbp]
\centering
\includegraphics[scale=0.5]{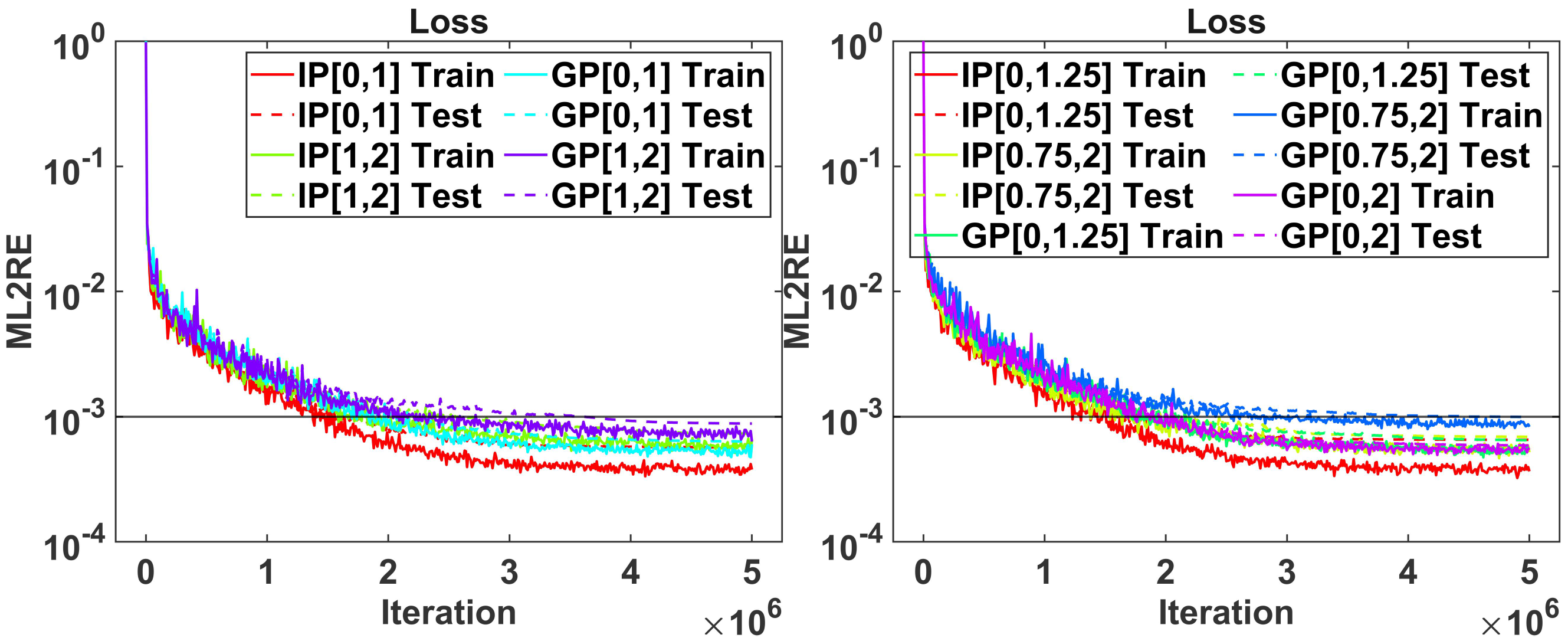}
\caption{Loss curve. (Left) D-R. (right) D-D. IP: interpolation. [0,2]: two-in-one.}
\label{fig_D-R-D-Loss}
\end{figure}

Figure \ref{fig_D-R-D-Loss} shows the loss curves of all subdomain networks under both schemes, with final losses below 1e-3. Error statistics are summarized in Table \ref{table_D-D-R}. Key observations:
\begin{itemize}
    \item Iterative convergence: DD-DeepONet outputs closely align with direct prediction results (when given subdomain BCs), confirming iterative convergence toward ground truth.
    \item Comparison of subdomain dataset acquisition methods: Interpolation-based method outperforms GP-based in accuracy. Precise interface data from interpolation enhances generalization, though neither method generates completely erroneous results. This depends on the generalization of the subdomain BCs that are supplemented during training.
    \item Error propagation: Suboptimal generalization contaminates neighboring solutions during iterations, degrading global accuracy.
    \item Data consolidation: Merging similar subdomains reduces the required quantity of subnetworks, while improving robustness through increased data complexity. More generally, sufficiently large subdomain datasets and adequately generalized subdomain networks can assemble PDE solutions over complex domains via DDM.
\end{itemize}

\renewcommand\tabularxcolumn[1]{>{\centering\arraybackslash}m{#1}}
\begin{table}[pos=htbp]
\centering
\caption{D-R and D-D scheme's results. DP: Direct prediction when given subdomain BC. 3D: DD-DeepONet. IP: interpolation. GP: Gaussian process.} 
\begin{tabularx}{\textwidth}{CCCCCC} 
\toprule
SCHEMES & DATA METHOD & DOMAIN & METHOD & ML2RE & MAE \\
\hline
    \multirow{8}{*}{  
        \begin{tabular}{@{}c@{}}
        Framework 1 \\ D-R \\ Non-Overlap
        \end{tabular} 
    } & \multirow{4}{*}{IP} & \multirow{2}{*}{$0 \leq y \leq 1$} & DP & 5.6632e-4 & 1.5878e-4  \\
    & & & 3D & 7.7601e-4 & 2.5248e-4  \\
    & & \multirow{2}{*}{$1 \leq y \leq 2$} & DP & 7.1894e-4 & 2.1155e-4  \\
    & & & 3D & 8.5320e-4 & 2.7349e-4  \\
\cline{2-6}
    & \multirow{4}{*}{GP} & \multirow{2}{*}{$0 \leq y \leq 1$} & DP & 1.4000e-3
 & 3.2993e-4  \\
    & & & 3D & 2.8000e-3 & 7.1346e-4  \\
    & & \multirow{2}{*}{$1 \leq y \leq 2$} & DP & 7.7779e-4 & 2.2160e-4  \\
    & & & 3D & 2.6000e-4 & 7.1277e-4  \\
\hline
    \multirow{12}{*}{  
        \begin{tabular}{@{}c@{}}
        Framework 1 \\ D-D \\ Overlap
        \end{tabular} 
    } & \multirow{4}{*}{IP} & \multirow{2}{*}{$0 \leq y \leq 1.25$} & DP & 6.6961e-4 & 1.6935e-4  \\
    & & & 3D & 9.3459e-4 & 3.0068e-4  \\
    & & \multirow{2}{*}{$0.75 \leq y \leq 2$} & DP & 6.7409e-4 & 2.0313e-4  \\
    & & & 3D & 7.8212e-4 & 2.5740e-4 \\
\cline{2-6}
    & \multirow{4}{*}{GP} & \multirow{2}{*}{$0 \leq y \leq 1.25$} & DP & 1.7000e-3 & 4.4338e-4  \\
    & & & 3D & 3.4000e-3 & 8.8248e-4  \\
    & & \multirow{2}{*}{$0.75 \leq y \leq 2$} & DP & 1.7000e-3 & 4.2230e-4  \\
    & & & 3D & 2.5000e-3 & 6.6252e-4  \\
\cline{2-6}
    & \multirow{4}{*}{
        \begin{tabular}{@{}c@{}}
            GP \\ Two-in-One
        \end{tabular} } & \multirow{2}{*}{$0 \leq y \leq 1.25$} & DP & 1.6000e-3 & 4.0489e-4  \\
    & & & 3D & 3.2000e-3 & 8.0687e-4  \\
    & & \multirow{2}{*}{$0.75 \leq y \leq 2$} & DP & 6.6289e-4 & 1.9119e-4  \\
    & & & 3D & 1.5000e-3 & 4.3130e-4  \\
\bottomrule
\end{tabularx}
\label{table_D-D-R}
\end{table}

\subsection{S2 examples}
In this section, we demonstrate learning the operator $\mathcal{G}: \mathbb{R}^{d} \rightarrow \Phi$ using two PDEs: 1. Laplace equation for a 3D resistance problem; 2. Incompressible Navier-Stokes equations for 2D steady pipe flow.

\subsubsection{Resistance}
With the advancement of VLSI technology, the feature size of integrated circuits (ICs) continues to shrink, and we have entered the era of nanometer processes, approaching the limits predicted by Moore's Law. At these process nodes, due to the continual reduction in the feature size of interconnects, they can no longer be simply treated as equipotential connections. It is necessary to account for the parasitic electromagnetic coupling effects between them \cite{2014Advanced}.

The extraction of parasitic resistance is crucial for IC design, as it significantly affects the performance of the circuit \cite{2014Advanced}. Excessive parasitic resistance can lead to signal delays and reduced speeds, increased power consumption, issues with signal integrity, and impaired circuit performance. This can also decrease the circuit's tolerance to noise, potentially cause local hotspots, and reduce the reliability and lifespan of the circuit.

With the advancement of technology, extraction techniques for 2D and 2.5D interconnects are increasingly unable to meet the demands, and the extraction of parasitic parameters for 3D interconnects has become progressively more important. However, due to the computational speed limitations of traditional numerical methods, industry-level parasitic extraction techniques do not rely entirely on traditional numerical methods for solving these problems. For example, for a straight conductor with two identical parallel ports, its resistance can be directly calculated using the analytical formula $R=\rho l / s = l / \sigma S$, where 
$l$ is the length along the current direction, $s$ is the cross-sectional area, and $\sigma$ is the conductivity. For more complex geometries, such as shapes with corners, traditional numerical methods may be required for calculation.

Parasitic resistance extraction requires solving a large number of electrostatic field equations with no free charges present, within a short time frame. For simplicity, the port where the current enters is set to a potential of 1, the port where the current exits is set to a potential of 0, and all other surfaces are considered insulated. The equations and BCs for all the problems to be solved are the same; the only difference is the changing geometry of the computation domain.

\begin{figure}[htbp]
\centering
\includegraphics[scale=0.5]{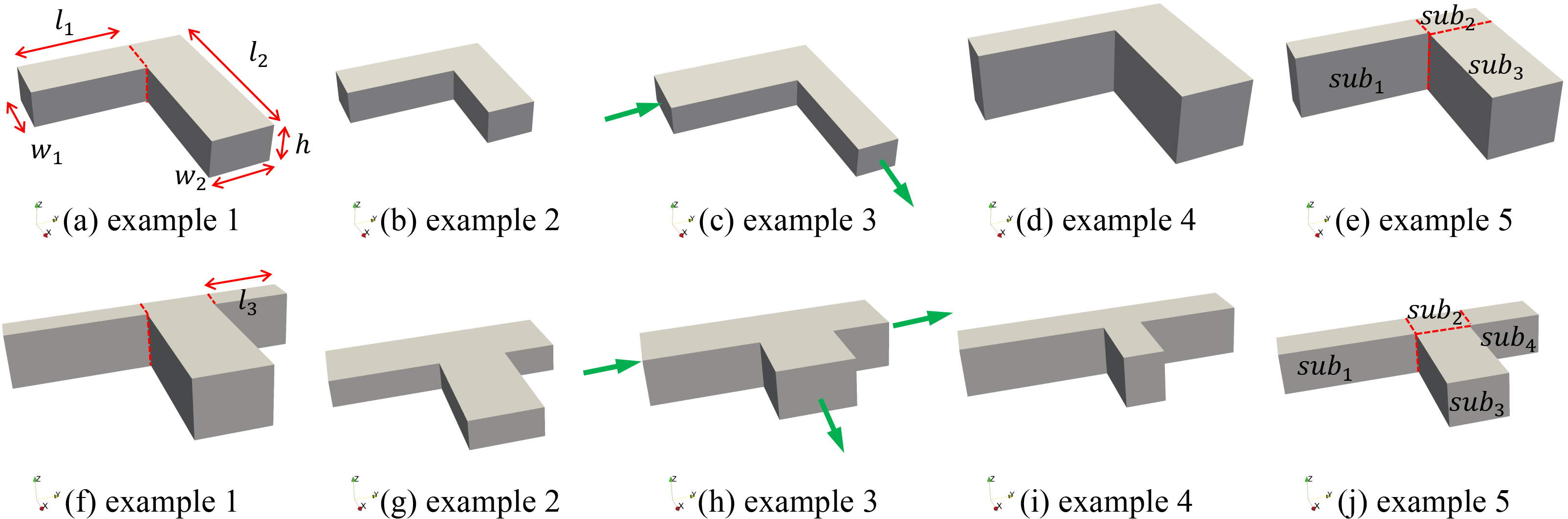}
\caption{L- and T-shaped computational domains with current ports marked by green arrows. Shape parameter ranges: $h, w_{1}, w_{2}$:1-3, $l_{1}$: 4-6, $l_{2}$: 5-8, $l_{3}$:2-5. The T-shape remains are same as the L-shape.}
\label{fig4.2_geo_LT}
\end{figure}

\begin{figure}[htbp]
\centering
\includegraphics[scale=0.65]{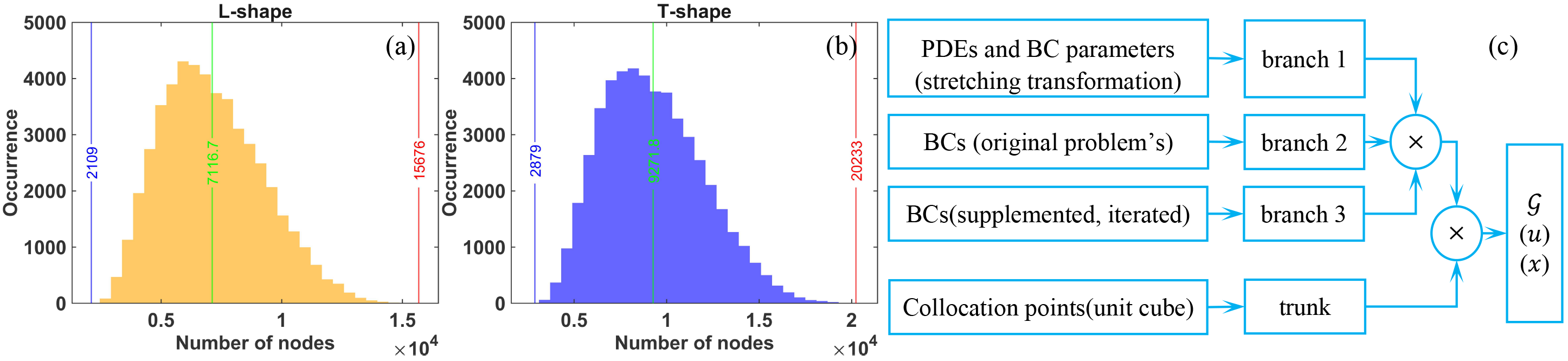}
\caption{(a, b)Histograms of the number of mesh nodes for L and T-shape. The vertical axis represents the frequency of sample occurrences. The red line is the maximum number of nodes, green is the average, and blue is the minimum. (c) Neural network structure and input.}
\label{fig4.2_geo_LT_node_ele}
\end{figure}

After obtaining the solution for the potential in the electrostatic field equation, the resistance $R_{IO}$ between the current entry and exit ports is calculated using Ohm's law:
\begin{equation}
    R_{IO}=\frac{V_{I}-V_{O}}{I_{O}}=\frac{1}{\int_{\Gamma_{O}} \sigma \frac{\partial u}{\partial \vec{n}} d \Gamma},
\label{eq4.2.1}
\end{equation}
$\Gamma_{O}$ is the current exit port, $\sigma = 5.998 \times 10^{7}$ is the conductivity, $V_{I}$ and $V_{O}$ is the electric potential, and $\vec{n}$ is the unit normal vector of the current exit ports boundary.

The ideal interconnect shape consists of structures formed by a series of cuboids. In this study, we focus on solving the resistances for L-shape and T-shape (Figure \ref{fig4.2_geo_LT}) interconnects in the same metal layer. For both shapes, we randomly generate 50,000 geometric samples. We use our custom-written $C++$ finite element code to generate data. All samples were discretized using second-order 10-node tetrahedral Lagrange elements. All numerical results were saved in VTK format files. The data at the subdomain interfaces are obtained using interpolation. Figure \ref{fig4.2_geo_LT_node_ele} shows the histogram of FEM grid node counts for all samples under both shapes.

\begin{figure}[htbp]
\centering
\includegraphics[scale=0.6]{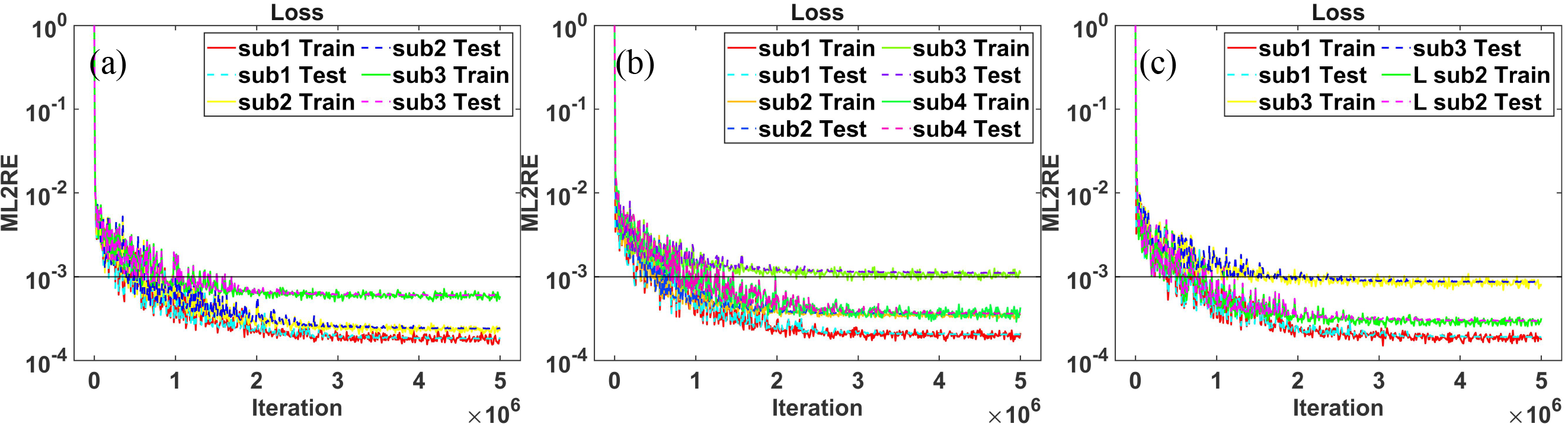}
\caption{Loss curve. (a) L-shape. (b) T-shape. (c) Loss curve after data merge.}
\label{fig_resistance}
\end{figure}

\begin{figure}[htbp]
\centering
\includegraphics[scale=0.55]{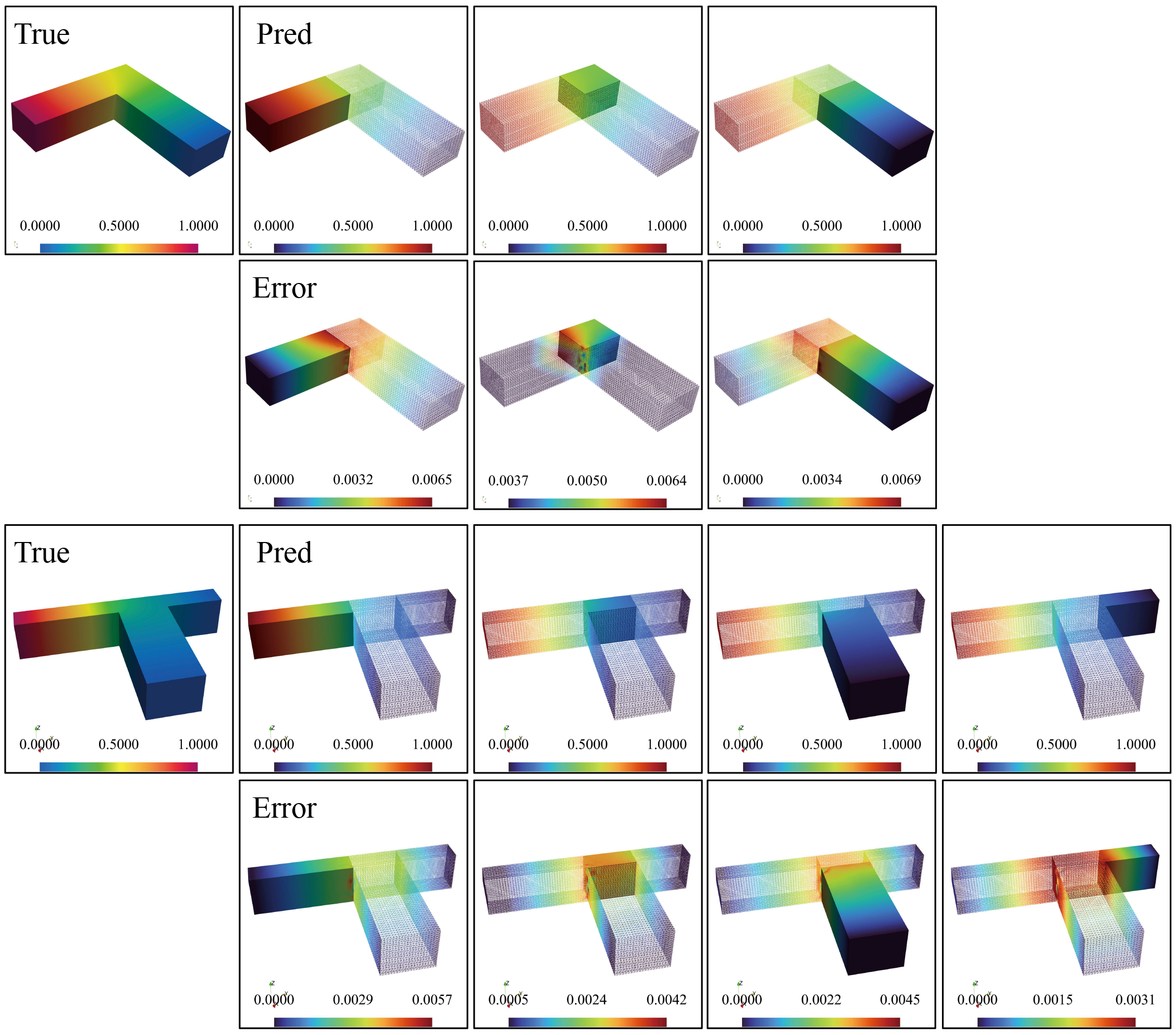}
\caption{The absolute error of an example in the L- and T-shape test set. The color range for both the true value and predicted value spans from the minimum to the maximum of the true value. For the error, the color range in each plot is scaled from its minimum to maximum value. }
\label{fig_LT_ae}
\end{figure}

We partition L- and T-shapes into 3 and 4 cuboid subdomains (Figure \ref{fig4.2_geo_LT}). Each subdomain is normalized to a unit cube via stretching transformations that embed shape parameters into PDEs coefficients and BCs. This converts shape-dependent problems into parameterized PDEs on fixed unit cubes, resolving DeepONet's input limitations on shape-dependent problems.

Each subdomain neural network adopts the architecture shown in Figure \ref{fig4.2_geo_LT_node_ele}(c). The number of product space inputs varies with subdomain configurations. We train neural networks for every L- and T-shape subdomain using their respective datasets. Additionally, Subdomains 1 and 3 exhibit similar configurations across both geometries. So we merge their datasets for joint training and retrain subdomain 2 in L-shape, while retaining trained networks for other subdomains.

After training the neural networks, we perform iterations in the test set using the coupling framework 2. To fully leverage GPU batch processing, we set 60 iterations for the L-shape and 50 for the T-shape. We set the relaxation factor $\theta = 0.5$, initial value on the interface equal to 0. The resistance value is obtained via numerical integration at the current outflow port.

Figure \ref{fig_resistance} shows the loss curves: Subdomain 3 in both geometries exhibits higher losses (stabilizing around 1e-3), while other subnetworks achieve lower values. Training and test losses remain closely aligned throughout. 

We evaluate results using: ML2RE and MAE of potential, mean relative error (MRE) of resistance, percentage of samples with resistance error (RE) <5\%:
\begin{equation}
    MRE = \frac{1}{N} \sum_{i}^{N} \frac{| R^{true}_i - R^{pred}_i |}{R^{true}_i}, RE = \frac{| R^{true} - R^{pred} |}{R^{true}}.
\end{equation}

\begin{table}[pos=htbp]
\centering
\caption{Statistical summary of resistance calculation results. Merge: the calculation results after data merging.} 
\begin{tabularx}{\textwidth}{CCCCCC} 
\toprule
 & MRE & RE$\le$5 & ML2RE & MAE & TIME \\
\hline
     L & 0.96 & 100$\%$ & 6.04e-3 & 2.99e-3 & 4.35$s$  \\
\cline{4-6}
     T front & 1.10 & 97.13$\%$ & \multirow{2}{*}{5.39e-3} & \multirow{2}{*}{1.83e-3} & \multirow{2}{*}{12.11$s$} \\
     T right & 1.18 & 96.50$\%$ & &  &  \\
\hline
     L (Merge) & 0.91 & 99.05$\%$ & 4.08e-3 & 1.91e-3 & 11.25$s$\\
\cline{4-6}
     T front (Merge) & 1.03 & 97.13$\%$ & \multirow{2}{*}{4.16e-3} & \multirow{2}{*}{1.34e-3}  & \multirow{2}{*}{13.98$s$}\\
     T right (Merge) & 1.11 & 99.82$\%$ &  &   &\\
\bottomrule
\end{tabularx}
\label{table_resistance}
\end{table}

Error statistics are shown in Table \ref{table_resistance}, where only the iteration time for potential calculations is reported. FEM computations took approximately 7077s for the L-shaped test case and 13378s for the T-shaped case, whereas DD-DeepONet takes less than 13 seconds. Figure \ref{fig_LT_ae} displays the potential fields computed via DD-DeepONet for L- and T-shapes, respectively.

\subsubsection{Pipe flow}

\begin{figure}[htbp]
\centering
\includegraphics[scale=0.6]{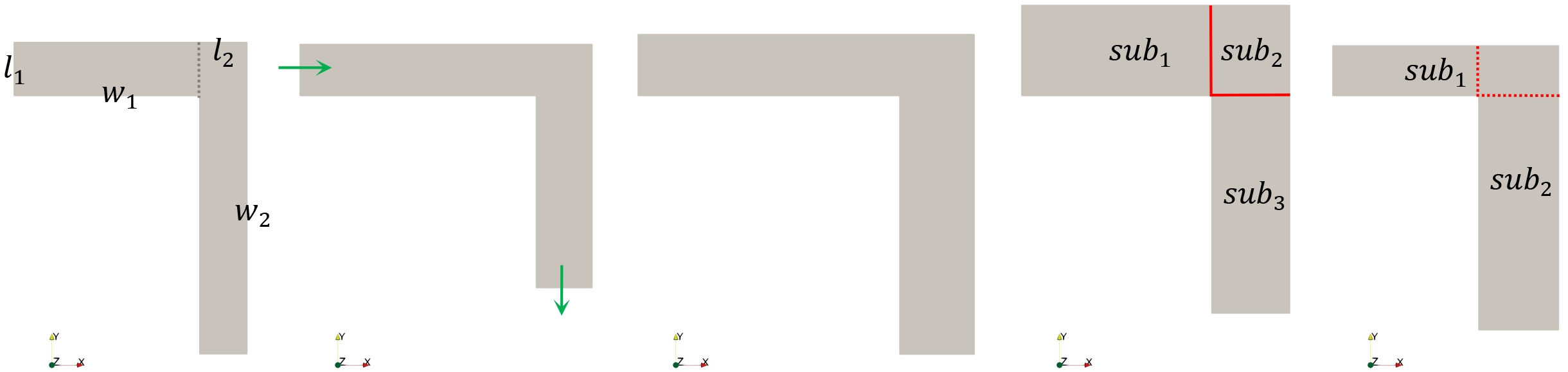}
\caption{L-shaped computational domains with flow direction of fluid by green arrows. Shape parameter ranges: $l_1, l_2$:1-3 $mm$, $w_{1}$: 4-8 $mm$, $w_{2}$: 6-10 $mm$.}
\label{Pipeflow_geo}
\end{figure}

We consider a 2D steady-state incompressible flow in an L-shaped pipe. The N-S equations are mathematically described as \cite{WOS:000517979600039}:
\begin{equation}
    \begin{aligned}
     \rho (\bm{u} \cdot \nabla) \bm{u} &= - \nabla p + \mu \nabla^{2} \bm{u}, \\
      \nabla \cdot \bm{u} &= 0. \\
    \end{aligned}
\label{eq_NS}
\end{equation}
where $p$ is the pressure, and $\bm{u}=(u,v)$ is the velocity vector. We set the fluid density $\rho = 1000 kg/m^3$ and dynamic viscosity $\mu = 0.001 Pa \cdot s$. At the inlet, the inflow velocity is set to $(0.1, 0) m/s$, while the outlet is configured as a fully developed flow. The remaining boundaries are set as no-slip walls. The geometric parameter ranges are shown in Figure \ref{Pipeflow_geo}. Using COMSOL, we generated 3428 samples with Reynolds numbers ranging from [100, 300], averaging 200.97. The scalar velocity $U = \sqrt{u^2+v^2}$ is the target for neural network learning.

\begin{figure}[htbp]
\centering
\includegraphics[scale=0.65]{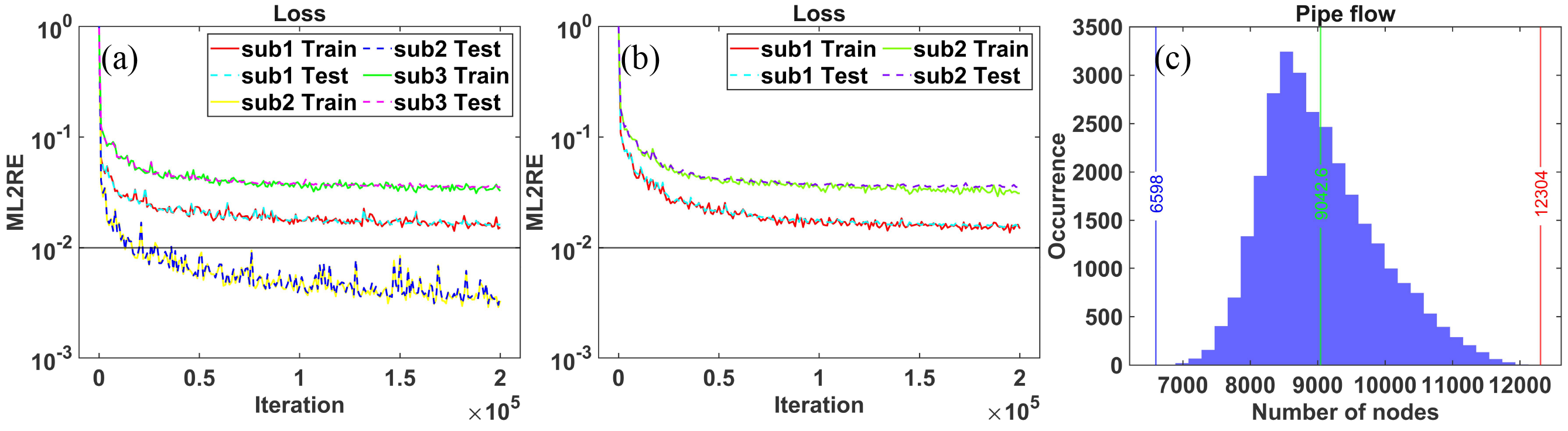}
\caption{Pipe flow results. (a) Loss curve of the coupling framework 2 (Non-overlap). (b) Loss curve of the coupling framework 1 (Overlap). (c) Histograms of the number of mesh nodes for pipe flow. The vertical axis represents the frequency of sample occurrences. The red line is the maximum number of nodes, green is the average, and blue is the minimum.}
\label{Pipeflow_loss}
\end{figure}

\begin{figure}[htbp]
\centering
\includegraphics[scale=0.5]{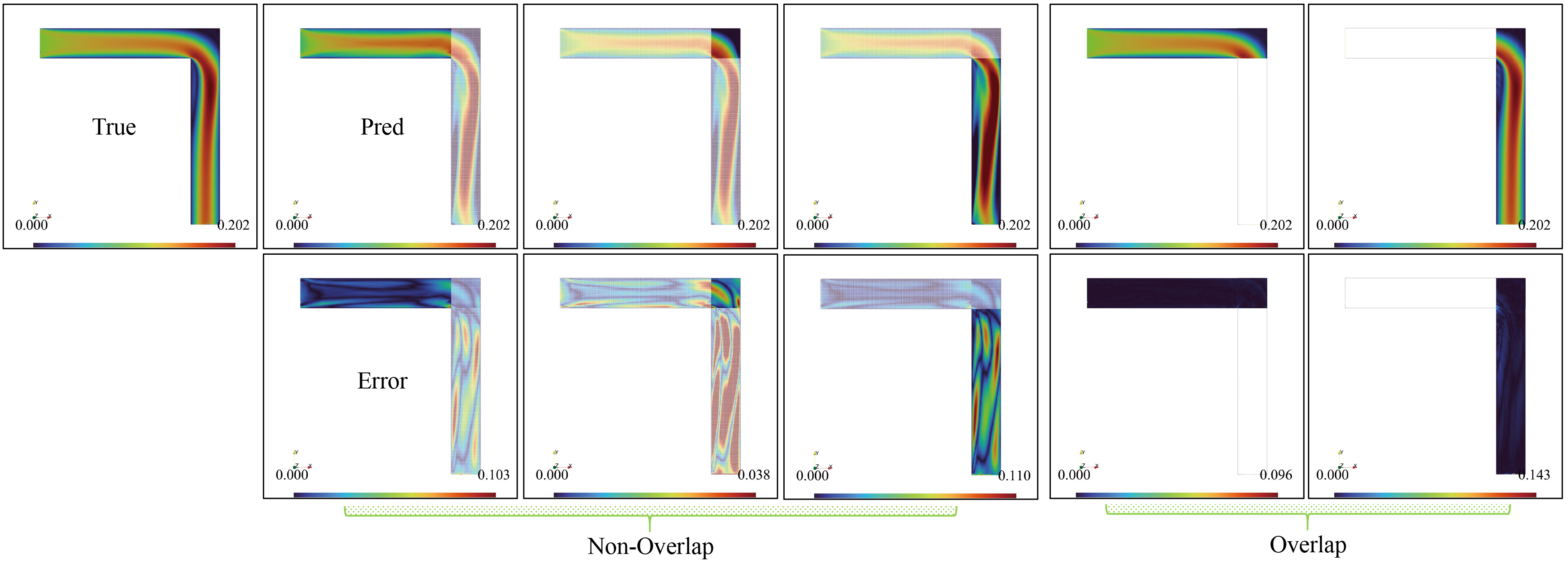}
\caption{The absolute error of an example in the test set. The color range for both the true value and predicted value spans from the minimum to the maximum of the true value. For the error, the color range in each plot is scaled from its minimum to maximum value.}
\label{Pipeflow_visual}
\end{figure}

We solve the problem using the non-overlapping coupling framework 2 and the overlapping D-D iterative scheme from the coupling framework 1. For the former, the computational domain is divided into three non-overlapping subdomains, while the latter uses two overlapping subdomains, as shown in Figure \ref{Pipeflow_geo}. For the non-overlapping case, the network input is similar to the resistance case. For the overlapping case, an additional proportion factor is included as input, representing the percentage of the supplemented BC relative to the subdomain boundary length. For both cases, we set 30 iterations, relaxation factor $\ theta=0.5$, and initial value on the interface equal to 0. 

Figures \ref{Pipeflow_loss} and \ref{Pipeflow_visual} present the computation results. It can be observed that the training and testing losses are higher in the subdomain near the outlet and lower near the inlet, while the non-overlapping subdomain 2 exhibits the lowest loss. This occurs because the flow changes at the bend, making the flow characteristics more complex and harder to learn compared to the inlet region. Non-overlapping subdomain 2, being smaller with a relatively uniform flow state, achieves the lowest loss.

On the test set, framework 2 achieves an ML2RE of 0.292 and an MAE of 0.025, while framework 1 achieves an ML2RE of 0.032 and an MAE of 0.0013. The two methods take 3.133 and 28.898 seconds of computation time, respectively. A comparison of the two iterative schemes shows that the overlapping method yields smaller errors. This is because framework 2 requires computing normal derivatives, which play a significant role in Robin BC and introduce larger errors. Furthermore, using Paraview to calculate the derivatives on non-grid nodes leads to an even greater increase in error. The overlapping iteration takes more computation time than the non-overlapping method, as it requires sampling interface points for each sample individually, preventing batch computation, unlike the non-overlapping approach. However, compared to the traditional numerical method's computation time of 3755s, the speedup is significant while maintaining good accuracy.

\subsection{S3 examples}
In this section, we demonstrate learning the operator $\mathcal{G}: \mathbb{R}^{d} \times \mathbb{R} \times \mathbb{R} \times \mathcal{P} \times \mathcal{Q} \rightarrow \Phi$ using two PDEs: 1. Drift-diffusion equation for a 2D PIN rib waveguide; 2. Multimedium Poisson equation for the interface problem.

\subsubsection{Drift-diffusion}

\begin{figure}[htbp]
\centering
\includegraphics[scale=0.55]{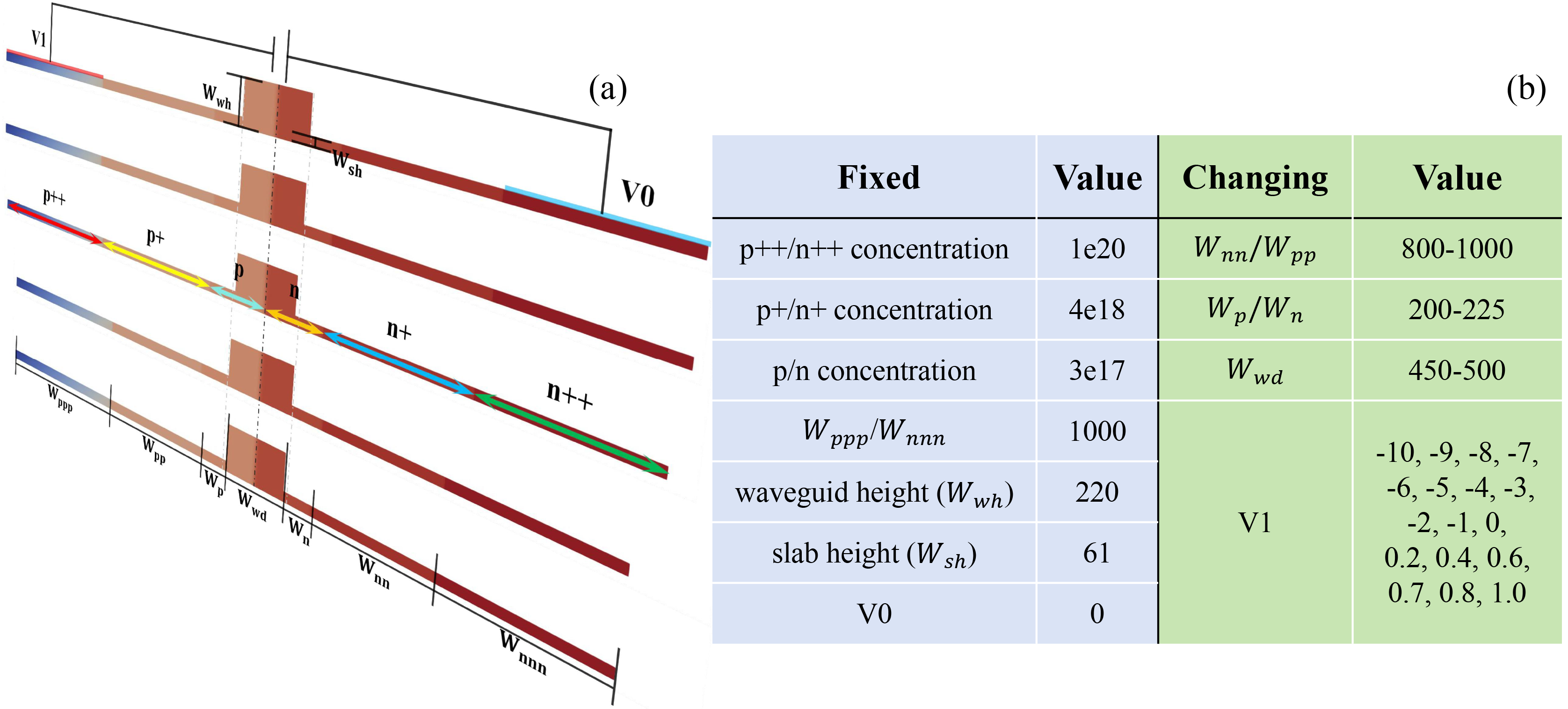}
\caption{(a) 2D view of the PIN rib waveguide. (b) Optimization variables and design parameters.}
\label{fig_DD_geo}
\end{figure}

The core physical processes in semiconductor silicon photonic devices (such as photodetectors, modulators, lasers) involve the generation of photogenerated carriers, their transport (drift and diffusion) under electric fields, recombination, and the resulting photocurrent or modulation of the optical field (refractive index). The drift-diffusion model is the most fundamental and widely used physical model for describing the transport behavior of these carriers (electrons and holes). 

Here, we consider a 2D steady-state drift-diffusion problem. The geometry and parameters are shown in Figure \ref{fig_DD_geo}, focusing only on the semiconductor domain \cite{WOS:001062861100006}. The structure consists of six regions with different P-type and N-type doping concentrations, each with distinct lengths. Bias voltage is applied uniformly from -10 to 1, sweeping from reverse to forward bias. The drift-diffusion equations are mathematically described as \cite{WOS:000793405100007}:
\begin{equation}
  \begin{aligned}
    \nabla \cdot (\epsilon \nabla \phi) &= -q (p - n + C), \\
    \nabla \cdot \bm{J_p} &= 0, \\
    \nabla \cdot \bm{J_n} &= 0, \\
    \bm{J_p} = -q n & \mu_n \nabla \phi + q D_n \nabla n, \\
    \bm{J_n} = -q p & \mu_p \nabla \phi - q D_p \nabla n, \\
    D_p = & \mu_p K_B T, \\
    D_n = & \mu_n K_B T, \\
  \end{aligned}
\label{eq_DD}
\end{equation}
where $\phi$ is the electrostatic potential, $\epsilon$ is the dielectric constant, $q$ is the fundamental electron charge, $p$ and $n$ are the electron
and hole concentrations inside the semiconductor. $C = N_D + N_A$ is the doping profile, which is assumed to be a given datum of the problem in terms of the donor and acceptor concentrations $N_D$(N-type doping) and $N_A$(P-type doping). $\bm{J_p}$ and $\bm{J_n}$ are hole and electron current density. The temperature $T$ of the crystal is constant. $K_B$ is the Boltzmann constant. $D_p$ and $D_n$ are hole and electron diffusion coefficients. $\mu_p$ and $\mu_n$ are hole and electron mobilities. 

\begin{figure}[htbp]
\centering
\includegraphics[scale=0.55]{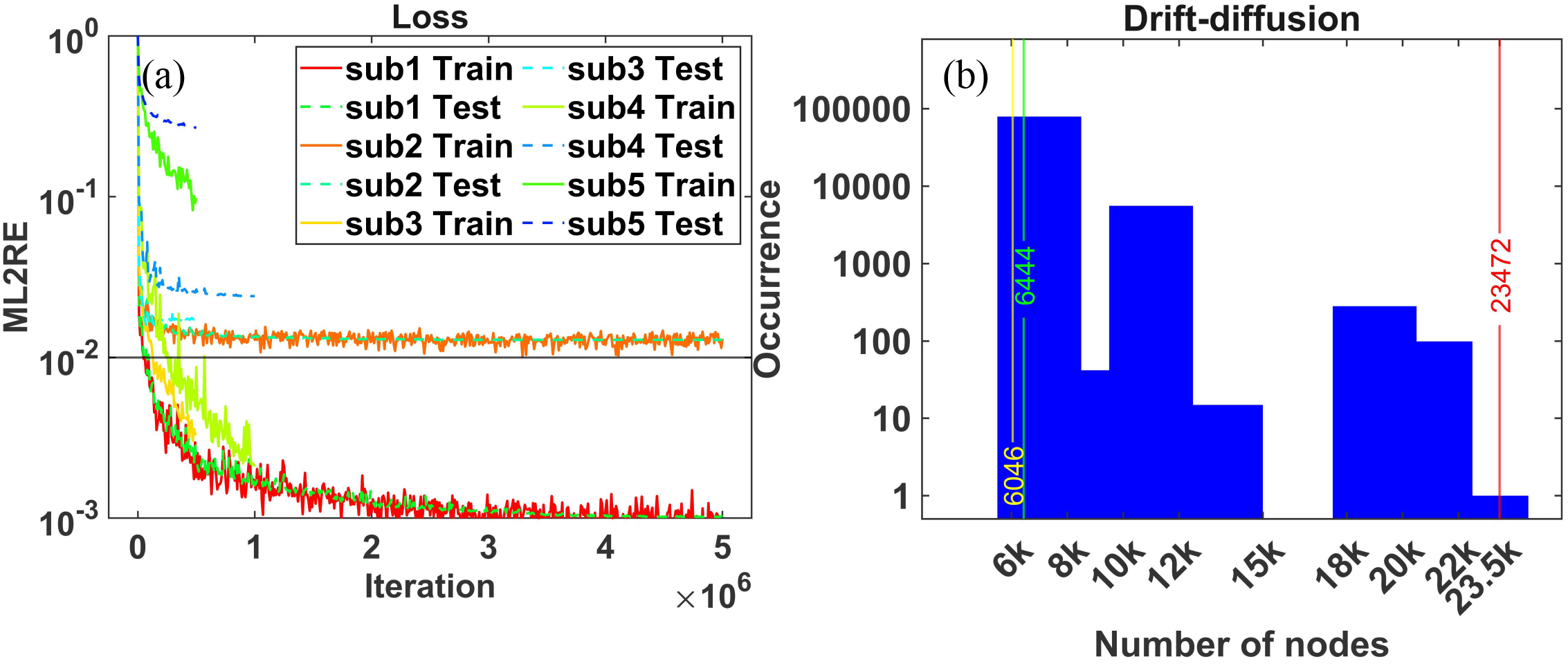}
\caption{(a) Loss curve of subnetwork for drift-diffusion. (b) Histograms of the number of mesh nodes for drift-diffusion. The vertical axis represents the frequency of sample occurrences. The red line is the maximum number of nodes, green is the average, and yellow is the minimum.}
\label{fig_DD_loss}
\end{figure}

\begin{figure}[htbp]
\centering
\includegraphics[scale=0.49]{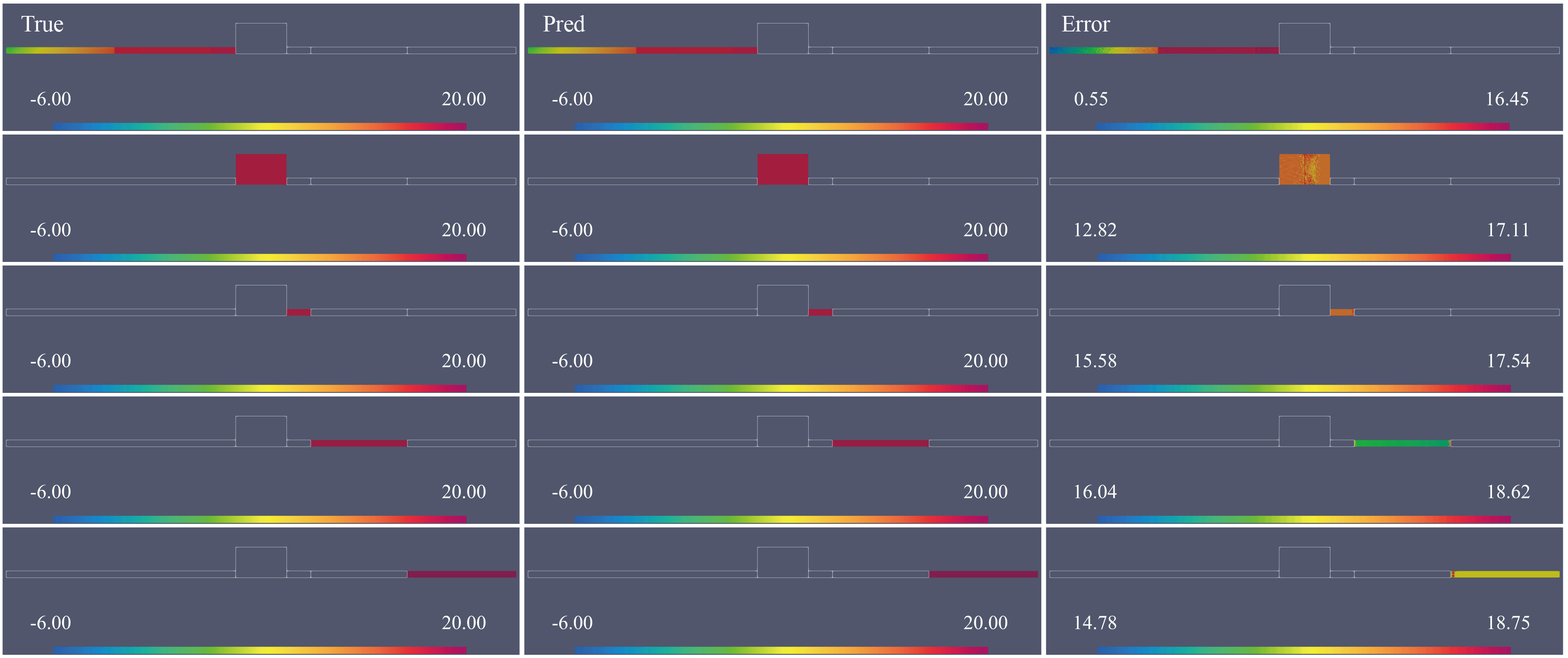}
\caption{The absolute error of a sample in the test set. The horizontal axis displays visualization results after applying $log_{10}$ to the data. The color range for both the true value and predicted value spans from the minimum to the maximum of the true value. For the error, the color range in each plot is scaled from its minimum to maximum value.}
\label{fig_DD_visual}
\end{figure}

We use COMSOL to generate data. The material within the computational domain is set to silicon, with default numerical settings in COMSOL. BC are set to insulation (homogeneous Neumann) and metal contacts (Non-homogeneous Dirichlet). We generated 5,000 geometric configurations, applying 17 different metal contact boundary voltages ($V_1$) to each, resulting in 85,000 samples. Grid node count statistics for all samples: see Figure \ref{fig_DD_loss} (b). We aim to learn the electron concentration $n$ by DD-DeepONet. Other variables are similar. 

We solve this problem using the non-overlapping coupling framework 2. The computational domain is divided into 5 subdomains as shown in Figure \ref{fig_DD_visual}. Given the data range $1 \times 10^{-6}$ to $1 \times 10^{20}$, we preprocess values near zero by taking $log_{10}(x+1)$ and larger values by taking $log_{10}(x)$, followed by z-score normalization. Subnetworks learn transformed training data. During testing, we apply the inverse transform to revert data for iterations, then re-transform to normalized data for predictions.

The neural network structure and inputs are similar to previous cases. Subdomain indices increase from left to right. Subdomains 1 and 2 were trained for 500k iterations, subdomains 3 and 5 for 50k iterations, and subdomain 4 for 100k iterations. The corresponding loss curves are shown in Figure \ref{fig_DD_loss}. The iteration count was set to 15, $\theta = 0.99$. The initial guess is taken from a sample in the training set. 

After 4s of iterative computation, the ML2RE reached 0.005235. Compared to COMSOL, more than a 3-day computation time, efficiency is significantly improved with good accuracy. Figure \ref{fig_DD_visual} shows a visualization result from the test set at 1V bias voltage, demonstrating satisfactory prediction performance.

\subsubsection{Multimedium Poisson}

\begin{figure}[htbp]
\centering
\includegraphics[scale=0.5]{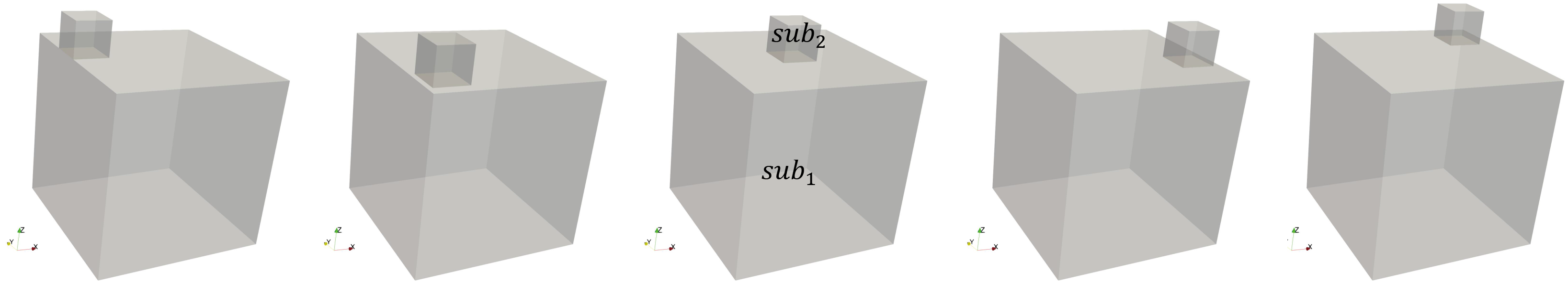}
\caption{Multimedium Poisson equation computational domain schematic.}
\label{fig_thermal_geo}
\end{figure}

We consider a Poission equation in two isotropic material layers, see Figure \ref{fig_thermal_geo}. The computational domain comprises two stacked cubes: cube 1 (lower, $l=1$), cube 2 (upper, $l=0.2$) with variable positioning. It describes material interface problems, applicable to: steady-state heat conduction, electrostatic fields, linear elastostatics, diffusion phenomena, among other practical problems. The mathematical description is as follows \cite{10247998}:

\begin{equation}
    \begin{aligned}
      -\epsilon_i \Delta u &= 1, \quad \text{in } \Omega, \, i=1,2 \\
      u &= 1, \quad \text{on } \Gamma_{\text{up}}, \\
      \epsilon_i \frac{\partial u}{\partial \vec{n}} &= 0, \quad \text{on } \Gamma_{\text{side}}, \\
      \epsilon_i \frac{\partial u}{\partial \vec{n}} &= 1 - u, \quad \text{on } \Gamma_{\text{down}}, \\
    \end{aligned}
\label{eq_thermal}
\end{equation}
with
\begin{equation}
    \begin{aligned}
        u_i &= u_j, \quad \text{on } \Gamma_{\text{interface}}, \\
        \epsilon_i \frac{\partial u_i}{\partial \vec{n}_i} &= -\epsilon_j \frac{\partial u_j}{\partial \vec{n}_j}, \quad \text{on } \Gamma_{\text{interface}}. \\
    \end{aligned}
\label{eq_thermal_with}
\end{equation}
$\epsilon_1 \in \{10.0, 10.1, 10.2, 10.3, 10.4, 10.5, 10.6, 10.7, 10.8, 10.9\}$ and $\epsilon_2 \in \{0.1, 0.2, 0.3, 0.4, 0.5, 0.6, 0.7, 0.8, 0.9, 1.0\}$  are the material conductivities. The two equations in (\ref{eq_thermal_with}) represent continuity conditions at the interfaces between different materials.

\begin{figure}[htbp]
\centering
\includegraphics[scale=0.65]{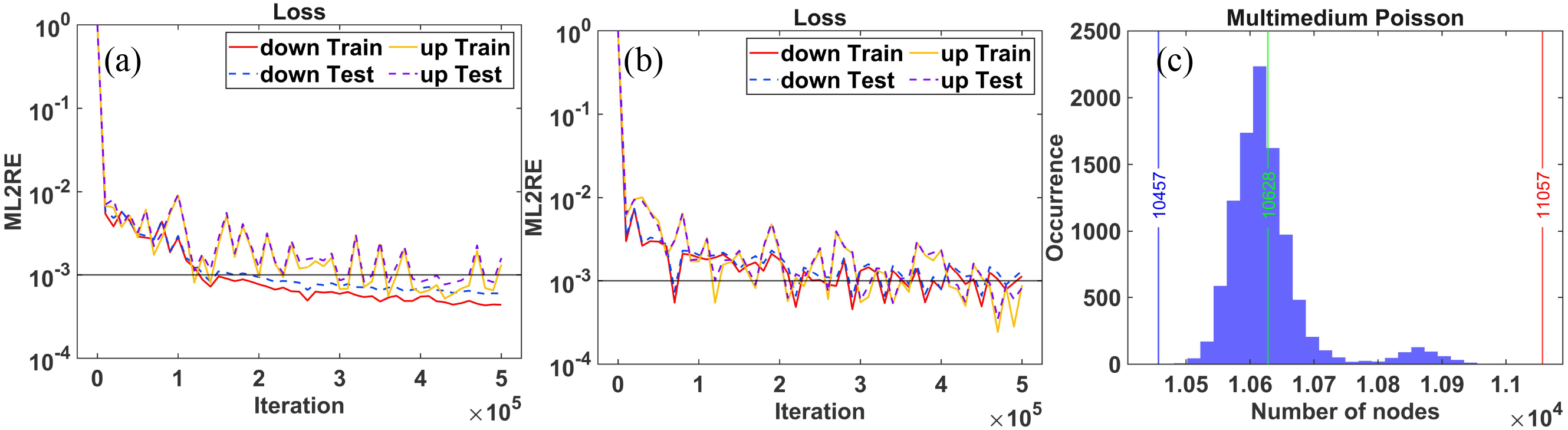}
\caption{Multimedium Poisson results. (a) Loss curve of the coupling framework 2. (b) Loss curve of iteration-free. (c) Histograms of the number of mesh nodes. The vertical axis represents the frequency of sample occurrences. The red line is the maximum number of nodes, green is the average, and blue is the minimum.}
\label{fig_thermal_loss}
\end{figure}

\begin{figure}[htbp]
\centering
\includegraphics[scale=0.5]{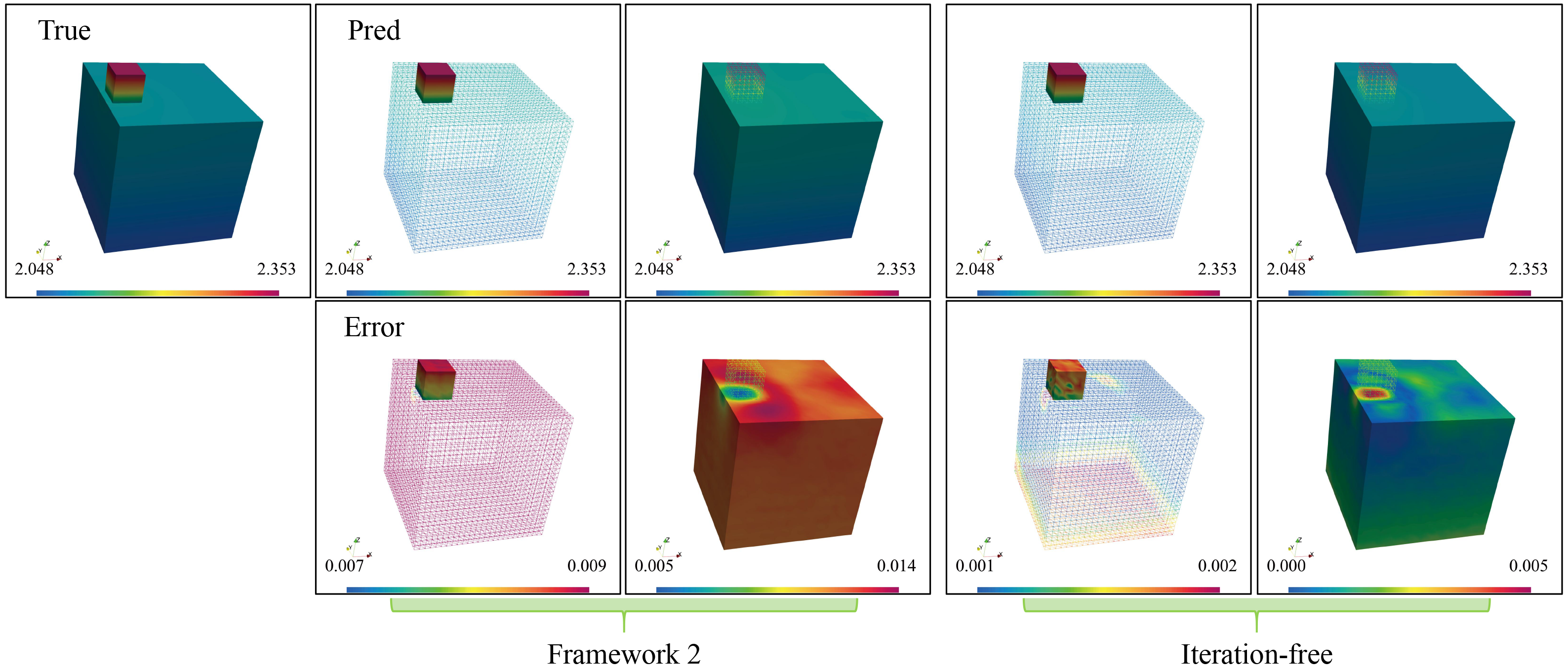}
\caption{The absolute error of an example in the test set. The color range for both the true value and predicted value spans from the minimum to the maximum of the true value. For the error, the color range in each plot is scaled from its minimum to maximum value.}
\label{fig_thermal_example}
\end{figure}

We generated 10,000 data samples using COMSOL. The problem was solved via the coupling framework 2 and iteration-free methods. Neural networks require additional input: cube 2's position on cube 1's surface. For cube 2's trunk input, we translate to $[0,0.2] \times [0,0.2] \times [0,0.2]$. For framework 2, we set the iteration count is 30, $\theta=0.5$, initial guess is all 2. For iteration-free, branch net input global domain data.

Figures \ref{fig_thermal_loss} and \ref{fig_thermal_example} show results from both methods. The training losses of both methods can reach 1e-3. On the test set, framework 2 achieves an ML2RE of 0.006024 and an MAE of 0.012726, while the iteration-free approach achieves an ML2RE of 0.001335 and an MAE of 0.002860. The iteration-free approach outperforms the coupling framework 2. Either because the iterative process extends beyond the training data and the model exhibits weaker extrapolation capability, or due to errors introduced by the Robin BC. However, both methods can achieve satisfactory accuracy overall.

\section{Conclusion, discussion and future work}

In this paper, we introduce the DD-DeepONet, which can simplify complex geometric PDEs assembled from cuboids/rectangles. We present its computation capabilities in three application scenarios that need repeated solutions in a short time. Through several PDE case studies, we evaluate several different DD-DeepONet iteration schemes. Subdomain network accuracy, initial guess, relaxation factor, and termination criteria significantly impact solution precision and computational time for the iteration scheme. For iteration-free, the influencing factor is only the subdomain network accuracy. We demonstrate that this method is practical and can alleviate the training difficulties. Combined with stretching transformations and translation, it can solve geometry-dependent problems and deliver solutions for linear/nonlinear PDEs quickly, which is remarkable for nonlinear PDEs. We highlight that both the DDM scheme and DeepONet in this method can be replaced to create entirely different structural frameworks from those presented in this paper, offering great scalability. 

For iteration schemes, they can serve as submodules to compose more complex PDE problems, unlike the iteration-free method. However, we note that interface values during iteration may fall outside the function space (training set), which can easily cause divergence. Therefore, to use them as submodules, sufficiently large and noise-resistant training sets are required, along with iteration schemes better suited to this approach, and neural networks with stronger generalization capability.

In the future, we will employ transfer learning to share knowledge across tasks with varying geometries, boundary conditions, or parameters, reducing the need for high-quality data and training costs. Physics-informed pretraining of DeepONet will be used to boost generalization and adaptability in complex settings. For large-scale simulations, a multiscale fidelity strategy can be applied—training coarse models for low-frequency structures, followed by fine-tuning with high-fidelity data. We also aim to extend the framework to multiphysics coupling, spatiotemporal PDEs, and stochastic PDEs, on complex geometries.



\section{Acknowledgements}
This work is supported in part by the National Key Research and Development Program of China (No. 2023YFB3001704), the Major Key Project of PCL (No. PCL2023A03), and the NSF of Fujian Province under Grants (No. 2024J09045). It is also supported in part by the Special Project for Research and Development in Key areas of Guangdong Province (No. 2021B0101190003) and the NSF of Guangdong Province (No. 2022A1515010831) through grants.

\clearpage
\appendix
\setcounter{section}{0}
\renewcommand{\thesection}{Appendix \Alph{section}}

\section{Neural network and dataset setup}
\label{Appendix_A}

\newcolumntype{C}{>{\centering\arraybackslash}X}
\begin{table}[pos=htbp]
\centering
\caption{Neural network and dataset setup of S1 example. SS: sample size. CPs: collocation points. DR: decay rate. Iter: iterations. IP: interpolation. 2[625, 256*4]: [[625, 256, 256, 256, 256],[625, 256, 256, 256, 256],[3, 256, 256, 256, 256]], other similiar. $DR(c)=0.999990 \times 10^{c}$. }
\begin{tabularx}{\textwidth}{CCCCCCCC} 
\toprule
\multicolumn{3}{c}{(train:test = 8:2)} & SS & CPs & NNs size & DR & Iter\\
\hline
\multicolumn{3}{c}{\multirow{6}{*}{Sample size}} & 5,000 & 30,625 & 2[625,256*4] & $DR(\frac{1}{8})$ & 5e6 \\
\cline{4-8}
 & & & 10,000 & \multirow{5}{*}{30,625} & \multirow{5}{*}{2[625,512*4]} & \multirow{5}{*}{$DR(\frac{1}{8})$} & \multirow{5}{*}{5e6}  \\
 & & & 20,000 &  &  &  &  \\
 & & & 30,000 &  &  &  &  \\
 & & & 40,000 &  &  &  &  \\
 & & & 50,000 &  &  &  &  \\
\hline
\multirow{2}{*}{Iteration-free} & \multicolumn{2}{c}{Non-DDM} & \multirow{2}{*}{50,000} & \multirow{2}{*}{99,937} & \multirow{2}{*}{2[1369,1024*4]} & \multirow{2}{*}{$DR(\frac{1}{6})$} & \multirow{2}{*}{3e6} \\
 & \multicolumn{2}{c}{DDM(Half)} & & & & &  \\
 \hline
\multirow{9}{*}{
        \begin{tabular}{@{}c@{}}
            Framework1 \\ Non-Overlap \\ and \\ Overlap
        \end{tabular} 
} & \multirow{2}{*}{D-R IP} & [0,1] & \multirow{2}{*}{30,000} & \multirow{2}{*}{15,625} & \multirow{2}{*}{2[625,512*4]} & \multirow{2}{*}{$DR(\frac{1}{8})$} & \multirow{2}{*}{5e6} \\
 & & [1,2] & & & & & \\
 \cline{2-8}
 & \multirow{2}{*}{D-R GP} & [0,1] & \multirow{2}{*}{30,000} & \multirow{2}{*}{15,625} & \multirow{2}{*}{2[625,512*4]} & \multirow{2}{*}{$DR(\frac{1}{8})$} & \multirow{2}{*}{5e6} \\
 & & [1,2] & & & & & \\
 \cline{2-8}
 & \multirow{2}{*}{D-D IP} & [0,1.25] & \multirow{2}{*}{30,000} & \multirow{2}{*}{19,375} & \multirow{2}{*}{2[625,512*4]} & \multirow{2}{*}{$DR(\frac{1}{8})$} & \multirow{2}{*}{5e6} \\
 & & [0.75,2] & & & & & \\
 \cline{2-8}
 & \multirow{3}{*}{D-D GP} & [0,1.25] & \multirow{2}{*}{30,000} & \multirow{3}{*}{19,375} & \multirow{3}{*}{2[625,512*4]} & \multirow{3}{*}{$DR(\frac{1}{8})$} & \multirow{3}{*}{5e6} \\
 & & [0.75,2] & & & & & \\
 & & Two-in-One & 60,000 & & & & \\
\bottomrule
\end{tabularx}
\label{table_NN_Setup_FEx}
\end{table}

\newcolumntype{C}{>{\centering\arraybackslash}X}
\begin{table}[pos=htbp]
\centering
\caption{Neural network and data set setup of resistance example. SS: sample size. CPs: collocation points. DR: decay rate. Iter: iterations. 2[441, 512*4]: [[9, 512],[441, 512, 512, 512, 512],[441, 512, 512, 512, 512],[3, 512, 512, 512, 512]], other similiar. $DR(c)=0.999990 \times 10^{c}$. }
\begin{tabularx}{\textwidth}{CCCCCCC} 
\toprule
\multicolumn{2}{c}{(train:test = 8:2)} & SS & CPs & NNs size & DR & Iter\\
\hline
\multirow{3}{*}{L-shape} & sub1 & \multirow{3}{*}{50,000} & 18081 & \multirow{3}{*}{2[441,512*4]} & \multirow{3}{*}{$DR(\frac{1}{8})$} & \multirow{3}{*}{5e6} \\
 & sub2 & & 9261 & & &  \\
 & sub3 & & 18081 & & &  \\
 \hline
\multirow{4}{*}{T-shape} & sub1 & \multirow{4}{*}{50,000} & 18081 & 2[441,1024*4] & \multirow{4}{*}{$DR(\frac{1}{8})$} & 5e6 \\
 & sub2 & & 9261 & 3[441,1024*4] & & 5e6  \\
 & sub3 & & 18081 & 2[441,1024*4] & & 4e6  \\
 & sub4 & & 18081 & 2[441,1024*4] & & 5e6  \\
 \hline
\multirow{3}{*}{Merge} & sub1 & \multirow{2}{*}{100,000} & 18081 & \multirow{3}{*}{2[441,1024*4]} & \multirow{3}{*}{$DR(\frac{1}{8})$} & \multirow{3}{*}{5e6} \\
 & sub3 & & 18081 & & &  \\
 & L-shape sub2 & 50,000 & 9261 & & &  \\
\bottomrule
\end{tabularx}
\label{table_NN_Setup_FEx_resistance}
\end{table}

\newcolumntype{C}{>{\centering\arraybackslash}X}
\begin{table}[pos=htbp]
\centering
\caption{Neural network and data set setup of pipe flow case. SS: sample size. CPs: collocation points. DR: decay rate. Iter: iterations. 2[40, 128*3]: [[6, 128],[40, 128, 128, 128],[40, 128, 128, 128],[2, 128, 128, 128]]; 2[40, 1, 512*3]: [[6, 512],[40, 512, 512, 512],[1, 512, 512, 512],[2, 512, 512, 512]], other similiar. $DR(c)=0.999990 \times 10^{c}$. }
\begin{tabularx}{\textwidth}{CCCCCCC} 
\toprule
\multicolumn{2}{c}{(train:test = 8:2)} & SS & CPs & NNs size & DR & Iter\\
\hline
\multirow{3}{*}{Non-overlap} & sub1 & \multirow{3}{*}{3428} & 120*40=4800 & \multirow{2}{*}{2[40, 128*3]} & \multirow{3}{*}{$DR(\frac{1}{8})$} & \multirow{3}{*}{2e5} \\
 & sub2 & & 40*40=1600 & & &  \\
 & sub3 & & 40*120=4800 & 1[40, 128*3] & &  \\
 \hline
\multirow{2}{*}{Overlap} & sub1 & \multirow{2}{*}{3428} & 160*40=6400 & \multirow{2}{*}{2[40, 1, 512*3]} & \multirow{2}{*}{$DR(\frac{1}{8})$} & \multirow{2}{*}{2e5} \\
 & sub2 & & 40*160=6400 & & &  \\
\bottomrule
\end{tabularx}
\label{table_NN_Setup_FEx_pipeflow}
\end{table}

\newcolumntype{C}{>{\centering\arraybackslash}X}
\begin{table}[pos=htbp]
\centering
\caption{Neural network and data set setup of drift-diffusion case. SS: sample size. CPs: collocation points. DR: decay rate. Iter: iterations. 2[9, 11, 2]: [[9, 512, 512, 512],[11, 512, 512, 512],[2, 512, 512, 512],[2, 512, 512, 512]], other similiar. $DR(c)=0.999990 \times 10^{c}$. }
\begin{tabularx}{\textwidth}{CCCCCCC} 
\toprule
\multicolumn{2}{c}{(train:test = 8:2)} & SS & CPs & NNs size & DR & Iter\\
\hline
\multirow{5}{*}{Non-overlap} & sub1 & \multirow{5}{*}{85000} & 201*11=2211 & 2[9, 11, 2] & \multirow{5}{*}{$DR(\frac{1}{8})$} & 5e6 \\
 & sub2 & & 101*61=6161 & 2[5, 11, 11, 2] & & 5e6 \\
 & sub3 & & 41*11=451 & \multirow{2}{*}{2[4, 11, 11, 2]} & & 5e5 \\
 & sub4 & & 101*11=1111 & & & 1e6 \\
 & sub5 & & 101*11=1111 & 2[4, 11, 2] & & 5e6 \\
\bottomrule
\end{tabularx}
\label{table_NN_Setup_FEx_DD}
\end{table}

\newcolumntype{C}{>{\centering\arraybackslash}X}
\begin{table}[pos=htbp]
\centering
\caption{Neural network and data set setup of multimedium Poisson case. SS: sample size. CPs: collocation points. DR: decay rate. Iter: iterations. $NN_1$: [[5, 256, 256, 256],[49, 256, 256, 256],[3, 256, 256, 256]]; $NN_2$: [[2, 256, 256, 256],[4, 256, 256, 256],[3, 256, 256, 256]]. $DR(c)=0.999990 \times 10^{c}$. }
\begin{tabularx}{\textwidth}{CCCCCCC} 
\toprule
\multicolumn{2}{c}{(train:test = 8:2)} & SS & CPs & NNs size & DR & Iter\\
\hline
\multirow{2}{*}{Framework 2} & sub1 & 10000 & 26*26*26=17576 & \multirow{2}{*}{$NN_1$}  & \multirow{2}{*}{$DR(\frac{1}{8})$} & \multirow{2}{*}{5e5} \\
 & sub2 & 10000 & 7*7*7=49 & & &  \\
 \hline
\multirow{2}{*}{Iteration-free} & sub1 & 10000 & 26*26*26=17576 & \multirow{2}{*}{$NN_2$} & \multirow{2}{*}{$DR(\frac{1}{8})$} & \multirow{2}{*}{5e5} \\
 & sub2 & 10000 & 7*7*7=49 & & &  \\
\bottomrule
\end{tabularx}
\label{table_NN_Setup_possion}
\end{table}

\clearpage
\bibliographystyle{unsrt}
\bibliography{cas-refs}



\end{document}